\documentclass[reqno]{amsart}
\usepackage{epsfig,latexsym,amsfonts,amssymb,amsmath,amscd,graphics,epic}
\usepackage{amsthm}
\usepackage[mathscr]{eucal}
\usepackage{mathrsfs}
\usepackage{mathbbol}
\usepackage{oldgerm,units}
\usepackage{wrapfig}
\usepackage{ifthen}

\newcommand{\T}{{\mathbb T}}
\newcommand{\Z}{{\mathbb Z}}\newcommand{\Q}{{\mathbb Q}}
\newcommand{\Inte}{{\operatorname{Int}}}\newcommand{\Ext}{{\operatorname{Ext}}}
\newcommand{\pr}{{\operatorname{pr}}}
\newcommand{\Del}{{\Delta}}
\newcommand{\const}{{\operatorname{const}}}\newcommand{\conv}{{\operatorname{conv}}}
\newcommand{\red}{{\operatorname{red}}}
\newcommand{\Real}{\mathbb R}\newcommand{\Vol}{{\operatorname{Vol}}}
\newcommand{\TrS}{\oplus}\newcommand{\Reg}{{\operatorname{Reg}}}
\newcommand{\TrP}{\odot}\newcommand{\bn}{\boldsymbol{n}}
\newcommand{\ba}{\boldsymbol{a}}\newcommand{\bw}{\boldsymbol{w}}
\newcommand{\bx}{\boldsymbol{x}}\newcommand{\bu}{\boldsymbol{u}}\newcommand{\bv}{\boldsymbol{v}}

\newcommand{\Net}{\mathbb N}
\newcommand{\eps}{\varepsilon}

\newcommand{\etype}[1]{\renewcommand{\labelenumi}{(#1{enumi})}}
\def\eroman{\etype{\roman}}

\def\wtU{\widetilde U}
\def\al{\alpha}
\def\dl{\delta}
\def\olK{\overline K}
\def\Sig{\Sigma}
\def\sig{\sigma}
\def\Lm{\Lambda}
\def\mcJ{{\mathcal J}}
\def\mcP{{\mathcal P}}

\def\v{x}
\def\u{u}
\def\a{A}
\def\b{B}
\def\c{C}
\def\var{\lambda}

\def\Var{ \Lm }

\newtheorem{theorem}{Theorem}[section]
\newtheorem{proposition}[theorem]{Proposition}
\newtheorem{definition}[theorem]{Definition}

\newtheorem{lemma}[theorem]{Lemma}
\newtheorem{conjecture}[theorem]{Conjecture}

\newtheorem{corollary}[theorem]{Corollary}

\newtheorem{remark}[theorem]{Remark}
\newtheorem{question}[theorem]{Question}

\begin{document}


\title[Idempotent Semigroups and Tropical Algebraic Sets] {Idempotent Semigroups and Tropical Algebraic Sets}


\author{Zur Izhakian}
\address{Department of Mathematics, Bar-Ilan University, 52900
Ramat-Gan, Israel} \email{zzur@math.biu.ac.il}

\author{Eugenii Shustin}
\address{School of Math. Sciences, Tel Aviv University, Ramat Aviv, 69978 Tel Aviv, Israel}
\email{shustin@post.tau.ac.il}






\begin{abstract}
The tropical semifield, i.e., the real numbers enhanced by the
operations of addition and maximum, serves as a base of tropical
mathematics. Addition is an abelian group operation, whereas the
maximum defines an idempotent semigroup structure. We address the
question of the geometry of idempotent semigroups, in particular,
tropical algebraic sets carrying the structure of a commutative
idempotent semigroup. We show that commutative idempotent
semigroups are contractible, that systems of tropical polynomials,
formed from univariate monomials, define subsemigroups with
respect to coordinate-wise tropical addition (maximum); and,
finally, we prove that the subsemigroups in $\Real^n$, which are
tropical hypersurfaces or tropical curves in the plane or in the
three-space, have the above polynomial description.
\end{abstract}

\maketitle

\section{Introduction}\label{intro}

Tropical geometry is a geometry over the tropical semifield
$\T=\Real\cup\{-\infty\}$ with the operations of tropical addition
and multiplication
$$a\TrS b=\max\{a,b\},\quad a\TrP b=a+b$$ (cf.
\cite{IMS,MikhalkinICM,MikhalkinBook,Sturmfels6366}). We equip
$\T^*=\Real$ with Euclidean topology, assuming that $\T$ is
homeomorphic to $[0,\infty)$. In this setting, tropical varieties
appear to be certain finite rational polyhedral complexes. The
simplest examples of tropical varieties, $\Real=\T^*$ and $\T$,
carry algebraic structures: for instance, $\Real$ is an abelian
group with respect to the tropical multiplication and is a
commutative idempotent semigroup with respect to tropical addition,
whereas $\T$ is a semigroup with respect to both the operations.
Thus, it is natural to ask about algebraic and geometric properties
of tropical varieties, equipped with one of these structures.

Group structure has addressed in
\cite{Gathmann:0601322,ZharkovMikhalkin}. In particular, the
tropical abelian varieties, {\it i.e.} tropical varieties which
are abelian groups whose operations are regular tropical
functions, say, tropical Jacobians, are just real tori (products
of circles and lines).

On the other hand, the tropical varieties (and more generally,
tropical algebraic sets) enhanced with a structure of an idempotent
semigroup, have not been touched yet. In this paper, we focus on the
geometric and algebraic properties of such tropical varieties. After
general consideration, resulting in the claim that connected
topological idempotent semigroups with a nontrivial center are
contractible (Theorem \ref{nt10}), we turn to an interesting
particular case of subsemigroups in $\Real^n$ equipped with the
coordinate-wise tropical addition. Observing that the tropical power
induces an endomorphism of $(\Real,\TrS)$, we conclude that tropical
polynomials consisting of only univariate monomials (termed
\emph{simple polynomials}) define subsemigroups of $(\Real^n,\TrS)$.
Yet, not every polyhedral complex which is a subsemigroup of
$(\Real^n,\TrS)$ can be defined by only simple polynomials. However,
we conjecture that such subsemigroups which are also tropical
varieties (called \emph{additive tropical varieties}) can be defined
by simple polynomials, and we prove this conjecture for the cases of
additive tropical hypersurfaces of arbitrary dimension (Theorem
\ref{nt1}) and for additive tropical curves in the plane and in the
three-space (Theorem \ref{thm:end}). As a consequence, we show that,
for any additive tropical variety, its skeletons support additive
tropical subvarieties (Theorem \ref{nt11}), and thus, all their
connected components are contractible.

\subsection*{Acknowledgements} The authors have been supported by
the Hermann-Minkowski Center for Geometry at the Tel Aviv
University, by a grant from the High Council for Scientific and
Technological Cooperation between France and Israel at the Ministry
of Science, Israel, and by the grant no. 448/09 from the Israeli
Science Foundation.

Part of this work was done during the authors' stay at the
Max-Planck-Institut f\"ur Mathematik (Bonn). The authors are very
grateful to MPI for the hospitality and excellent working
conditions.

\section{Topology of idempotent semigroups}

A topological semigroup is a pair $(U,\psi)$, where $U$ is a
topological space, and $\psi:U\times U\to U$ is continuous and
associative, i.e.,
$$\psi(u,\psi(v,w))=\psi(\psi(u,v),w),\quad u,v,w\in U\ .$$ In the sequel, we consider
only topological semigroups and therefore will omit the word
``topological" and write semigroups, for short. Moreover, when the
operation is clear from the context, we will write $U$ for
$(U,\psi)$. Also, we will often write $uv$ instead of $\psi(u,v)$;
no confusion will arise.

The \emph{center} of a semigroup $(U,\psi)$ is defined to be the
set
$$C(U,\psi):=\{u\in U\ :\ \psi(u,v)=\psi(v,u)\ \text{for all}\
v\in U\}\ .$$ A semigroup $(U,\psi)$ is called \emph{idempotent}
if $\psi(u,u)=u$ for all $u\in U$.

We start with two simple observations.

\begin{lemma}\label{nl1}
Any connected component of an idempotent semigroup forms  a
subsemigroup.
\end{lemma}

\begin{proof}
Let $(U,\psi)$ be an idempotent semigroup with $U_0\subset U$ a
connected component. Then, $\psi(U_0\times U_0)$ is connected,
and,  since the diagonal remains in $U_0$, it is  contained in
$U_0$.
\end{proof}

\begin{lemma}\label{nl2}
Any commutative idempotent semigroup $U$ is a directed poset with
respect to the relation
\begin{equation}v\succ u\quad\Longleftrightarrow\quad v=au\ \text{for some}\
a\in U\ ,\label{ne3}\end{equation} which in its turn is compatible
with the semigroup operation in the following sense:
$$v\succ u\ \Longrightarrow\ vw\succ uw\ \text{for
all}\ w\in U\ .$$
\end{lemma}

\begin{proof}
Reflexivity and transitivity of relation (\ref{ne3}) are
immediate. Next, if $u\succ v$ and $v\succ u$, then $u=va$,
$v=ub$, and we obtain $u=va=uab=uabb=ub=v$. Hence, relation
(\ref{ne3}) defines a partial order. Since $u\prec uv$ and $v\prec
uv$ for any $u,v\in U$, we obtain a directed set. Finally,
$$v\succ u\ \Longrightarrow\ v=au\ \Longrightarrow vw=a(uw)\
\Longrightarrow\ vw\succ uw\ .$$
\end{proof}

Our main observation is the following:

\begin{theorem}\label{nt10} Let $(U,\psi)$ be an idempotent
semigroup with a nonempty center, and let $U$ be a connected
topological space homotopy equivalent to a CW-complex. Then, $U$
is contractible.
\end{theorem}

\begin{proof}
By assumption, there exists $u_0\in C(U,\psi)$. We shall show that
$\pi_k(U,u_0)=0$ for all $k\ge 1$. This will yield the
contractibility by the classical Whitehead theorem.

Represent elements of $\pi_k(U,u_0)$ by maps $\gamma:I^k\to U$,
where  $I=[0,1]$, $\gamma(\partial I^k)=u_0$, taken up to homotopy
relative to $\partial I^k$. The operation of $\pi_k(U,u_0)$ is
then induced by the map composition of $\gamma_1,\gamma_2:I^k\to
U$ with $\gamma_1(\partial I^k)=\gamma_2(\partial I^k)=u_0$,
defined as
$$\gamma_1*\gamma_2:I^k\to U,\quad\gamma_1*\gamma_2(t_1,\dots,t_k)=\begin{cases}
\gamma_1(2t_1,t_2,\dots,t_k),\quad&0\le t_1\le\frac{1}{2},\\
\gamma_2(2t_1-1,t_2,\dots,t_k),\quad&\frac{1}{2}\le t_1\le 1\ .
\end{cases}$$
For each map $\gamma:I^k\to U$ with $\gamma(\partial I^k)=u_0$, we
have:
$$\gamma=\psi(\gamma,\gamma)\stackrel{h}{\sim}\psi(u_0*\gamma,\gamma*u_0)=
\psi(u_0,\gamma)*\psi(\gamma,u_0)=\psi(u_0,\gamma)*\psi(u_0,\gamma)\
,
$$ where  $\stackrel{h}{\sim}$ stands for the homotopy relative to
$\partial I^k$, and $u_0$ means the constant map. Furthermore,
$$\psi(u_0,\gamma)\stackrel{h}{\sim}\psi(u_0,\psi(u_0,\gamma)*\psi(u_0,\gamma))=
\psi(u_0,\psi(u_0,\gamma))*\psi(u_0,\psi(u_0,\gamma))=\psi(u_0,\gamma)*\psi(u_0,\gamma)\
; $$ which altogether means that $\gamma\stackrel{h}{\sim}u_0$, and
we are done.
\end{proof}

\begin{remark}\label{nr1} The hypothesis on the non-emptiness of the center cannot be
discarded  from Theorem \ref{nt10}, since, for example, in any
topological space $U$ one can define the structure of an
idempotent semigroup by letting $\psi(u,v)=u$ for any $u,v\in U$.
\end{remark}

A natural question arising from the preceding discussion is:

\medskip
\noindent {\bf Question:} {\it Does any contractible topological
space admit a structure of an idempotent semigroup with a nonempty
center?}

\medskip

This is indeed so in the following particular situation.

\begin{proposition}\label{np1}
Any 1-dimensional contractible CW-complex admits a structure of a
commutative idempotent semigroup.
\end{proposition}

\begin{proof}
Let $U$ be a 1-dimensional contractible CW-complex and pick a point
$u_0\in U$. For any point $u\in U$ there is a unique path
$\gamma_u\subset U$ joining $u$ with $u_0$ that is homeomorphic
either to $I=[0,1]$, or to a point, according as $u\ne u_0$ or
$u=u_0$. The intersection of two paths $\gamma_u$ and $\gamma _v$,
$u,v\in U$, is a path $\gamma_w$, for some $w\in U$, and thus,  we
define
$$uv=w\quad\text{as far as}\quad \gamma_u\cap\gamma_v=\gamma_w\
.$$ It is easy to check that this operation is associative,
commutative, and idempotent.
\end{proof}

Contractible tropical curves are called \emph{rational}
\cite{MikhalkinEnumerative}. Accordingly, our results show that a
tropical curve carries a structure of a commutative idempotent
semigroup iff it is rational. This contrasts with the case of
compact tropical curves having the structure of an abelian group:
these are elliptic (or more precisely, are homeomorphic to a
circle \cite{ZharkovMikhalkin}).

\section{Basic tropical algebraic geometry}

For the reader's convenience we first recall some necessary
definitions and facts in tropical geometry; these can be found in
\cite{AR,Einsiedler8311,Gathmann:0601322,GKM,IMS,Mikhalkin04,Sturmfels6366}.
We also introduce some new notions to be used throughout the text.

\subsection{Tropical polynomials and tropical algebraic sets}
A tropical polynomial is an expression of the form
\begin{equation}\label{ne0}
 f=
\bigoplus_{\omega\in\Omega}\a_\omega\TrP
\var_1^{\omega_1}\TrP\cdots\TrP \var_n^{\omega_n},
\end{equation} where
$\Omega\subset\Z^n$ is a finite nonempty set of points
$\omega=(\omega_1,\dots,\omega_n)$ with nonnegative coordinates,
$\a_\omega\in\Real$ for all $\omega\in\Omega$; here and in the
sequel, the power $a^m$ means $a$ repeated $m$ times, i.e.,
$a^m=\underbrace{a\TrP\cdots\TrP a}_m = ma$. We write a polynomial
as $f= \bigoplus_{\omega\in\Omega}\a_\omega\TrP \Var^{\omega}$,
where $\Var^{\omega}$ stands for
\mbox{$\var_1^{\omega_1}\TrP\cdots\TrP \var_n^{\omega_n}$}, and
denote the semiring of tropical polynomials by $\T[\Var]$. Abusing
notation, we will sometimes write $f(\var_{i_1},\dots,
\var_{i_m})$ for $f \in \T[\Var]$, indicating that $f$ involves
only the variables $\var_{i_1},\dots, \var_{i_m}$.

Any tropical polynomial $f\in \T[\Var]\backslash\{-\infty\}$
determines a piecewise linear convex function
\mbox{$f:\Real^n \to \Real$}:
\begin{equation}f =\bigoplus_{\omega\in\Omega}\a_\omega\TrP
\Var^\omega\ \longmapsto\
f(\bu)=\max_{\omega\in\Omega}(\langle\bu,\omega\rangle+A_\omega) \ ,
\label{ne1}\end{equation} where $\bu$ stands for the $n$-tuple
$(u_1,\dots, u_n) \in \Real^n$, and $\langle*,*\rangle$ is the
standard scalar product. Unlike the classical polynomials over an
infinite field, here the map of $\T[\Var]$ to the space of functions
is not injective. Some of the linear functions in the right-hand
side of \eqref{ne1} can be omitted without changing the function; we
call the corresponding monomials of $f$ \emph{inessential}, while
the other monomials are called \emph{essential}.

\begin{remark}
We denote a tropical polynomial and the corresponding function by
the same symbol, no confusion will arise. Whenever we write an
expression with formal variables $\lambda_i$, we assume a
polynomial, otherwise we mean a function. The value of the
function, corresponding to a polynomial $f\in\T[ \Lm ]$ at a point
$\bu\in\T^n$ is denoted by $f(\bu)$ or $f( \Lm )\big|_{\bu}$.
\end{remark}

Given a tropical polynomial $f$,  $Z_{\T}(f)$ is defined to be the
set of points $\bu\in\T^n$ on which the value $f(\bu)$ is either
equal to $-\infty$, or attained by at least two of the monomials
on the left-hand side of (\ref{ne1}). When
$f\in\T[\Var]\backslash\{-\infty\}$ is nonconstant, the set
$Z_\T(f)$ is a proper nonempty subset of $\T^n$, and is called
\emph{an affine tropical hypersurface}. Note that for the constant
polynomial $f= -\infty$ we have $Z_\T(-\infty)=\T^n$.

Letting $I=\langle f_1,\dots,f_s\rangle\subset\T[\Var]$ be a
finitely generated ideal, the set
$$Z_\T(I):=\bigcap_{f\in
I}Z_\T(f)\subset\T^n$$ is called \emph{an affine tropical
(algebraic) set}. Clearly, $Z_\T(I)=Z_\T(f_1)\cap\dots\cap
Z_\T(f_s)$. Indeed, taking a polynomial $f=g_1\TrP f_1\oplus \dots
\oplus g_n\TrP f_n\in I$ and a point $\bu\in
Z_\T(f_1)\cap\dots\cap Z_\T(f_s)$, we have $ f(\bu)= g_i(\bu)+
f_i(\bu)$ for some $i=1, \dots ,n$. Suppose $ f(\bu)\ne -\infty$.
Then, the value of $ g_i(\bu)$ is attained by a monomial $\b_\mu
\TrP \Var^\mu$ of $g_i$ and the value of $ f_i(\bu)$ is attained
by some pair of monomials of $f_i$, say $\a_{\omega'}\TrP
\Var^{\omega'}$ and $\a_{\omega''}\TrP \Var^{\omega''}$. Thus, the
value of $ f(\bu)$ is attained by the two monomials $\a_{\omega'}
\TrP \b_\mu \TrP \Var^{\omega' + \mu}$ and $\a_{\omega''} \TrP
\b_\mu \TrP \Var^{\omega''+ \mu}$ of $f$; that is $\bu \in
Z_\T(f)$.

It is more convenient (and traditional) to consider tropical
algebraic sets in $\Real^n\subset\T^n$ (a tropical torus, cf.
\cite{MikhalkinBook}). So, for a tropical polynomial
$f\in\T[\Var]\backslash\{-\infty\}$, we let
$$Z(f):=Z_\T(f)\cap\Real^n.$$ This set can be viewed as the corner
locus of the function $f$, i.e., the set of points $\bu\in\Real^n$
on which $f$ is not differentiable, or, equivalently, the set of
points $\bu\in\Real^n$ where the value $f(\bu)$ is attained by at
least two of the linear functions in the right-hand side of
(\ref{ne1}). For example, $Z(f)$ is nonempty as long as $f$
contains at least two monomials.
Given a finitely generated ideal $I=\langle
f_1,\dots,f_s\rangle\subset\T[x_1,\dots,x_n]$, the set
$Z(I):=Z_\T(I)\cap\Real^n$ is called a \emph{tropical (algebraic)
set in $\Real^n$}.

\subsection{Tropical varieties}\label{nsec3} A \emph{finite polyhedral complex} (briefly, \emph{FPC})
in $\Real^n$ is a pair $(P, \mcP )$, where $P\subset\Real^n$ and $
\mcP $ is a finite set of distinct convex closed polyhedra in
$\Real^n$, such that:
\begin{itemize}\item $P=\bigcup_{ \sig \in \mcP } \sig $;
\item
if $ \sig \in \mcP $, then any proper face of $ \sig $ also
belongs to $ \mcP $;
\item if $ \dl , \sig \in \mcP $, then $ \dl \cap \sig $ is
either empty, or is a common face (not necessarily proper) of $
\dl $ and $ \sig $.
\end{itemize} Let $\dim(P, \mcP )=\max\{\dim( \sig ):  \sig \in \mcP \}$.
An FPC $(P, \mcP )$ is said to be \emph{pure-dimensional} if any $
\dl \in \mcP $ is a face of some $ \sig \in \mcP $ with $\dim(
\sig )=\dim(P, \mcP )$. An FPC $(P, \mcP )$ is called
\emph{rational}, if all the linear spaces
$$\Real \sig  : =\{\bu-\bu'\ :\ \bu,\bu'\in \sig \},\quad \sig \in \mcP \ ,$$ are defined over
$\Q$.

It is not difficult to see that tropical sets are rational FPC and
vice versa.

An  $m$-dimensional \emph{tropical variety} in $\Real^n$, $n>m$,
is a rational FPC $(P, \mcP )$ of pure dimension $m$ equipped with
the weight function $w$ which is defined on the set of
top-dimensional cells of $(P, \mcP )$, gives positive integral
values, and satisfies the balancing condition at any cell $\tau\in
\mcP $ of dimension $m-1$:
\begin{equation}\sum w( \sig )\bv_\tau( \sig )=0\in\Z^n/\Z\tau\
,\label{ne2}\end{equation} where the sum is taken over all
$m$-dimensional $ \sig \in \mcP $ containing $\tau$ as a face,
$\Z\tau=\Real\tau\cap\Z^n$, and $\bv_\tau( \sig )$ is a generator
of the lattice $\Z \sig /\Z\tau$ oriented to the cone centered at
$\tau$ and directed by $ \sig $.

In this paper we deal mainly with a weaker notion of tropical
variety which we call a \emph{tropical set-variety}:

\begin{definition}
Let $(P, \mcP ,w)$ be a tropical variety in $\Real^n$. We call the
set $P$ a tropical set-variety.
\end{definition}

Namely, when working with tropical set-variety,  we get rid of the
weight function and the FPC-structure. However, a tropical
set-variety can be canonically represented as the union of convex
polyhedra. Given an $m$-dimensional tropical set-variety $P$, we
denote by $\Reg(P)$ the set of points of $P$ where $P$ is locally
homeomorphic to $\Real^m$.

\begin{lemma}\label{ln1}
Let $P$ be an $m$-dimensional tropical set-variety. Then,
\begin{itemize}\item the closures of the connected components of $\Reg(P)$ are rational
$m$-dimensional convex polyhedra; \item if $K_1,K_2$ are two
connected components of $\Reg(P)$, and $\dim(\olK_1\cap
\olK_2)=m-1$, then $ \sig = \olK_1\cap \olK_2$ is a common face of
$ \olK_1$ and $ \olK_2$.\end{itemize}
\end{lemma}

\begin{proof} The case of $m=1$ is evident, and we assume that $m\ge 2$.

Suppose that the closure $ \olK$ of a connected component $K$ of
$\Reg(X)$ is not convex, that is there are closed convex polyhedra
$ \sig ,\tau,\xi\subset\partial \olK$, $\dim( \sig
)=\dim(\tau)=m-1$, $\dim(\xi)=m-2$, $\xi= \sig \cap\tau$, such
that  $ \olK$ is not convex in a neighborhood of a point
$\bx\in\Inte(\xi)$. Without loss of generality, we may assume that
$P$ is a cone with vertex $\bx$ (the weight function and the
balancing condition will be naturally inherited by the cone from
any structure of tropical variety on $P$).

Take an $(n-m+2)$-dimensional subspace $V$ of $\Real^n$ defined
over $\Q$, passing through $\bx$ and transverse to $\xi$. It
supports a tropical variety with one cell of weight $1$ whose
intersection with $P$ is a two-dimensional tropical set-variety
(see \cite{AR,GKM} for details) which possesses a connected
component $K\cap V$ of its regular part with a non-convex closure
$ \olK\cap V$; more precisely, this component is the complement of
the convex sector $S$ spanned by the rays $ \sig \cap V$ and
$\tau\cap V$ in the two-plane $\Pi=\bx+\Real K\cap V$.

 Let $W\subset\Real^n$ be a hyperplane defined over $\mathbb Q$
(again a tropical set-variety) containing the plane $\Pi$ and
transverse to each one-dimensional cell of $P\cap V$ which is not
parallel to $\Pi$. Then $P\cap V\cap(a+W)$, for a small generic
vector $a\in\Real^n$, is a tropical set-curve, whose projection to
$\Pi$ (being a plane tropical set-curve, push-forward in
terminology of \cite{AR,GKM}) is contained in a neighborhood of
the sector $S$, which  contradicts the balancing condition.

Now suppose that, for some two connected components $K_1,K_2$ of
$\Reg(P)$, $ \sig = \olK_1\cap \olK_2$ is not a common face and
has dimension $m-1$. Two situations are possible: either $ \sig
\subset\partial \olK_1\cap\partial \olK_2$, or $ \sig \cap\Inte(
\olK_1)\ne\emptyset$. They both can be viewed as limit case of the
preceding consideration when either $ \sig ,\tau$ lie in the same
$(m-1)$-face of $ \olK_1$, or $ \sig =\tau$. Then, the above
argument literally leads either to a plane tropical set-curve
different from a straight line and lying in a half-plane, or to a
plane set-curve in a neighborhood of a ray. Both the cases
contradict the balancing condition.
\end{proof}

\begin{definition}\label{d34}
For an $m$-dimensional tropical set-variety $X$ (in $\Real^n$),
denote by $X^{(m-1)}$ the union of the $(m-1)$-dimensional faces
of the closures of the connected components of $\Reg(X)$.
\end{definition}

\begin{question}\label{q1} Is $X^{(m-1)}$ a tropical set variety?
\end{question}

The answer is yes for tropical set-curves (evident), tropical
set-hypersurfaces (commented in the next section), and for additive
tropical set-varieties as we show in section \ref{sec-last}.

\subsection{Tropical hypersurfaces}\label{nsec2}
An important example of a tropical variety is a \emph{tropical
hypersurface}, i.e. a tropical variety in $\Real^n$ of dimension
$n-1$. By \cite[Proposition 2.4, Corollary 2.5]{Mikhalkin04}, for
any tropical hypersurface $(P, \mcP ,w)$ in $\Real^n$, there
exists a tropical polynomial $f(\lambda_1, \dots ,\lambda_n)$
which satisfies $P=Z(f)$ and possesses a number of properties
listed below.

If $f$ is given by \eqref{ne0}, then the Legendre dual to $f$
function $\nu_f$ is defined on the Newton polytope $\Delta_f$ of
$f$ (the convex hull of the set $\Omega$ in \eqref{ne0}), is
convex and piece-wise linear. The graph of $\nu_f$ can be viewed
as the lower part of the convex hull of the set
$\{(\omega,-\a_\omega)\in\Real^{n+1}\ :\
\omega\in\Delta\cap\Z^n\}$, whre $\a_\omega$ being the
coefficients from formula (\ref{ne0})). The maximal linearity
domains of $\nu_f$ and their faces (which all are convex lattice
polytopes) define an FPC-structure $S(f)$ on $\Delta_f$. This
structure is dual to the FPC-structure $ \Sig _P$ on $\Real^n$,
given by $ \mcP $ and the closures of the connected components of
$\Real^n\backslash P$. Namely, there is a one-to-one
correspondence between the cells of $S(f)$ and $ \Sig _P$ which
inverts the incidence relation and is such that
\begin{itemize}
\item the vertices of $S(f)$ on $\partial\Delta_f$ correspond to the
closures of the unbounded components of $\Real^n\backslash P$, and
the vertices of $S(f)$ in $\Inte(\Delta_f)$ correspond to the
closures of the bounded components of $\Real^n\backslash P$, \item
the cells of dimension $m>0$ in $S(f)$ correspond to cells of
dimension $n-m$ in $ \mcP $, and the corresponding cells are
orthogonal, \item the weight of an $(n-1)$-dimensional cell of ${
\mcP }$ equals the lattice length of the dual segment of $S(f)$.
\end{itemize}

It follows immediately that the vertices of the subdivision
$S({f})$, or, equivalently, the components of $\Real^n\backslash P$,
correspond bijectively to the essential monomials of $f$; in
particular, the vertices of $\Delta$ always correspond to the
essential monomials of $f$. Another immediate consequence is that
the Newton polytopes of tropical polynomials $g$ such that $Z(g)=P$
and their FPC structure $S(g)$ have the same combinatorial type so
that the corresponding cells are parallel.

In connection to Question \ref{q1}, we recall the following
well-known fact, supplying it with a simple proof.

\begin{lemma}
The proper faces of the closures of the connected components of
the complement in $\Real^n$ to a tropical set-hypersurface $P$
define a FPC structure $ \mcP $ on $P$, and, for any $k=0, \dots
,n-2$, the set $P^{(k)}=\bigcup_{ \sig \in \mcP ,\dim( \sig ) \le
k} \sig $ is a $k$-dimensional topical set-variety.
\end{lemma}

\begin{proof}
If $k=0$, then $P^{(0)}$ is a finite set which is always a
zero-dimensional tropical set-variety. So, fix $0<k<n-1$ and Let
$P=Z(f)$ for some tropical polynomial $f$. Let $w_k$ be the weight
function $w_k( \sig )=\Vol_\Z( \sig ^*)$, $ \sig \in \mcP $,
$\dim( \sig )=k$, where $ \sig ^*$ is the dual polytope in the
subdivision $S(f)$ of the Newton polytope $\Delta$ of $f$, and
$\Vol_\Z( \sig ^*)$ is the lattice volume of $ \sig ^*$ (i.e. the
ratio of the Euclidean $k$-dimensional volume $\Vol_\Real( \sig
^*)$ and $\Vol_\Real(\Delta_ \sig )$, $\Delta_ \sig $ being the
minimal lattice simplex in $\Real \sig $).

 We will show that
$P^{(k)}$ with the FPC structure $ \mcP ^{(k)}=\{ \sig \in{ \mcP
}\ :\ \dim( \sig )\le k\}$ and the weight function $w_k( \sig )$
is a $k$-dimensional topical set-variety. Hence, we pick $\tau\in{
\mcP }$, $\dim (\tau) =k-1$, and prove that
$$\sum_{\renewcommand{\arraystretch}{0.6}
\begin{array}{cc}
\scriptstyle{ \sig \in \mcP ,\ \dim ( \sig )=k}\\
\scriptstyle{\tau\subset \sig }
\end{array}}\Vol_\Z( \sig ^*)\cdot\bv_\tau( \sig )=0\in\Z^n/\Z\tau\ ,$$ or,
equivalently,
\begin{equation}\sum_{\renewcommand{\arraystretch}{0.6}
\begin{array}{cc}
\scriptstyle{ \sig \in \mcP ,\ \dim( \sig ) =k}\\
\scriptstyle{\tau\subset \sig }
\end{array}}\Vol_\Z( \sig ^*)\cdot\bv^{\perp}_\tau( \sig )=0\ ,\label{eh4}\end{equation} where
$$\bv_\tau( \sig )=\bv^{\perp}_\tau( \sig )+\bv^{\parallel}_\tau( \sig ),\quad\bv^{\parallel}_\tau( \sig )
\in\Real\tau,\ \bv^{\perp}_\tau( \sig )\perp\Real\tau\ .$$ Notice
that
$$\bv^{\perp}_\tau( \sig )=\frac{\Vol_\Real(\Delta_ \sig )}{k\Vol_\Real(\Delta_\tau)}\bn_\tau( \sig )\
,$$ where $\bn_\tau( \sig )$ is the unit vector in $\Real \sig $
orthogonal to $\Real\tau$ and directed inside $ \sig $. Observing
that $\Vol_\Real(\Delta_ \sig )=\Vol_\Real(\Delta_{ \sig ^*})$, we
rewrite (\ref{eh4}) as
\begin{equation}\sum_{\renewcommand{\arraystretch}{0.6}
\begin{array}{cc}\scriptstyle{\tau\subset \sig \in \mcP }
\\ \scriptstyle{\dim ( \sig ) =k}
\end{array}}\Vol_\Real( \sig ^*)\cdot\bn_\tau( \sig )=0\ .\label{eh5}\end{equation}
Finally, observe that $\bn_\tau( \sig )$ is the outer normal in
$\Real \tau^*$ to the facet $ \sig ^*$ of the polytope $\tau^*$,
and hence (\ref{eh5}) turns into the polytopal Stokes formula.
\end{proof}

\section{Simple Additive Tropical Sets}\label{nsec1}

Subsemigroups in $(\Real^n,\TrS)$ which are tropical algebraic
sets are called \emph{additive tropical sets}.

\begin{lemma}\label{nl3}
Let $u_1, \dots, u_n \in \Real$, then $\left(\bigoplus_{i=1}^n u_i
\right)^s = \bigoplus_{i=1}^n u_i^s$ for all $n,s \in \Net$.
\end{lemma}

\begin{proof} We apply the double induction on $s$ and $n$. Fixing
$n=2$, the case $s=1$ is evident. Then, the induction step from
$s-1$ to $s$ (where $s\ge 2$) goes as follows:
$$
\begin{array}{lll}
  (u_1 \TrS u_2)^s & = & (u_1 \TrS u_2)\TrP(u_1 \TrS u_2)^{s-1}
\\[1mm]
 & =& (u_1
\TrS u_2)\TrP(u_1^{s-1} \TrS u_2^{s-1}) \\[1mm]
& = & u_1^s \TrS u_1^{s-1}\TrP u_2 \TrS u_1\TrP u_2^{s-1} \TrS
u_2^{s}\ .
\end{array}
 $$
When $u_1 = u_2$  the required equality is clear: $(u_1\TrS
u_1)^s=u_1^s=u_1^s\TrS u_1^s$. If $u_1
> u_2$, then $u_1^s
> u_1^{s-1}\TrP u_2 \TrS u_1\TrP u_2^{s-1} \TrS u_2^{s}$, and hence $(u_1\TrS u_2)^s  = u_1^s=u_1^s\TrS u_2^s$.
The case of $u_1<u_2$ is treated similarly. The proof is then
completed by an induction on $n$.
\end{proof}

A tropical polynomial $f\in\T[\Var]$ is called \emph{simple} if
each of its monomials is univariate, or a constant.

\begin{corollary}\label{thm:powerOfPoly}
Any tropical algebraic set $Z(I)\subset\Real^n$, where the ideal
$I\subset\T[\Var]$ is finitely generated by simple polynomials, is
additive.
\end{corollary}

\begin{proof}
It is sufficient to prove that, for any simple polynomial
$f\in\T[\Var]$, the set $Z(f)\subset\Real^n$ is closed under the
operation $\TrS$.

Given $\bu,\bv\in Z(f)$, by Lemma \ref{nl3}, we have $f(\bu\TrS
\bv)=f(\bu)\TrS f(\bv)=\max\{f(\bu),f(\bv)\}$. Suppose that
$f(\bu\TrS \bv)=f(\bu)$. Since $\bu\in Z(f)$, then
$f(\bu)=M_1(\bu)=M_2(\bu)$, for some two distinct monomials $M_1$
and $M_2$ of $f$. By our assumption and by Lemma \ref{nl3}, we
have
$$f(\bu\TrS\bv)=f(\bu)=M_1(\bu)=M_2(\bu)\ge\max\{M_1(\bv),M_2(\bv)\}\ ,$$ and
hence $M_i(\bu\TrS \bv)=M_i(\bu)\TrS M_i(\bv)=M_i(\bu)=f(\bu\TrS
\bv)$, $i=1,2$; that is $\bu\TrS \bv\in Z(f)$.
\end{proof}

An additive tropical set of the form $Z(I)$ with an ideal
$I\subset\T[\Var]$ finitely generated by simple polynomials, is
called a \emph{simple additive tropical set}.

Not all additive tropical sets are simple; for example, the
horizontal ray $$R=\{(t,0)\ :\ t\ge 0\}\subset\Real^2$$ is a
tropical algebraic set defined by the ideal $$I=\langle \var_1\TrP
\var_2\TrS \var_1\TrS \var_2,\ \var_1\TrP \var_2\TrS \var_1\TrS
(-1)\TrP \var_2\rangle\subset\T[\var_1,\var_2]\ ,$$ and it is
additive. On the other hand, $R$ is not simple. Indeed, due to the
duality, described in Section \ref{nsec2},  the Newton polygon of
a tropical polynomial $f \in\T[\var_1,\var_2] $ such that
$Z(f)\supset R$ must have a (vertical) side with the outer normal
$(1,0)$. For a simple polynomial $f$ with two variables, which may
contain only monomials of the form $\a_0$, $\a_i\TrP \var_1^i$, or
$\b_j\TrP \var_2^j$, this is possible only when $f=f(\var_2)$,
i.e., $\var_1$ is not involved in $f$. But then, $Z(f)$ must
contain the whole straight line through the ray $R$, and so does
$Z(I)$ for $I$ generated by such simple polynomials.

However, we propose the following converse to Corollary
\ref{thm:powerOfPoly}.

\begin{conjecture}\label{nc1}
Any additive tropical set-variety in $\Real^n$ is simple.
\end{conjecture}

Next we prove this conjecture for the three particular cases:
tropical set-hypersurfaces, affine subspaces of $\Real^n$, and
tropical set-curves in $\Real^2$ and $\Real^3$.

\section{Additive tropical set-hypersurfaces and affine subspaces}

\begin{theorem}\label{nt1}
A tropical set-hypersurface $P\subset\Real^n$ is additive if and
only if\; $P=Z(f)$ for some simple tropical polynomial
$f\in\T[\Var]$.
\end{theorem}

\begin{proof} It is sufficient to prove the ``only if" implication.

{\it Step 1}. Let $P=Z(f)$ be additive. Without loss of generality,
since multiplication by a monomial and removal of inessential
monomials does not affect $Z(f)$ (see details in Section
\ref{nsec2}), we may assume that $f$ is not divisible by any
monomial and  it contains only essential monomials. Then, all the
monomials of $f$ are encoded by points lying   on the boundary of
the Newton polytope $\Delta$ of $f$. Indeed, otherwise we would have
an essential monomial, corresponding to a vertex of the subdivision
$S({f})$ in $\Inte(\Delta)$, and thus it is dual to a bounded
component of $\Real^n\backslash Z(f)$. But, in view of Theorem
\ref{nt10},  the latter is impossible, since $Z(f)$ is contractible
while the boundary of a bounded component of $\Real^n\backslash P$
would give a nontrivial $(n-1)$-cycle in $Z(f)$.

{\it Step 2}. Since $P\ne\emptyset$, $f$ has at least two monomials.
Let $\a_\omega\TrP \Var^\omega$ and $\a_\tau\TrP \Var^\tau$ be two
monomials of $f$ such that $\omega\ne\tau$, where
$\omega,\tau\in\Z^n$,
 and the corresponding hyperplane, cf. \eqref{ne1},
\begin{equation}\langle \bu,\omega\rangle+\a_\omega=\langle
\bu,\tau\rangle+\a_\tau\label{ne4}\end{equation} contains an
$(n-1)$-dimensional cell $D$ of $P$. We claim that the $n$-tuple
$\omega-\tau$ has at most two nonzero coordinates, and the product
of any pair of coordinates of $\omega-\tau$ is non-positive. Indeed,
otherwise, one could write the equation describing (\ref{ne4}) as
$a_1 \var_1+\dots+a_n \var_n=b$ with $a_i,a_j>0$ for some $i\ne j$,
and then one could choose two sufficiently close points
$\bu'=(u'_1,\dots,u'_n)$, $\bu''=(u''_1,\dots,u''_m)$ in the
interior of $D$ for which
$$u'_i>u''_i,\quad u'_j<u''_j,\quad u'_k=u''_k,\ \text{for all}\
k\ne i,j\ .$$ But then, $\bu'\TrS \bu''\not\in D$, since this sum
does not satisfy (\ref{ne4}).

{\it Step 3}. Suppose that $n=2$, and $P=Z(f)$ for
$f\in\T[\var_1,\var_2]$.

Let $f$ contain the monomials $\a_i\TrP \var_1^i$ and $\b_j\TrP
\var_2^j$ with $i,j>0$. Assuming that $f$ has a (essential)
monomial $\a_{kl}\TrP \var_1^k\TrP \var_2^l$ with some $k,l>0$,
and taking into account the conclusions of Step 1, we obtain the
three vertices $(i,0)$, $(0,j)$, and $(k,l)$ of the subdivision
$S({f})$ lying on the boundary of the Newton polygon $\Delta$.
Along the conclusion of Step 2, the sides of $\Delta$ cannot be
directed by vectors with positive coordinates, and hence the
tropical curve $U$ necessarily has either
\begin{itemize}\item  a pair of rays, directed by vectors
with negative coordinates (see Figure \ref{nfig1}(a),(b)), or
\smallskip
\item  a pair of rays, directed by vectors with positive
coordinates (see Figure \ref{nfig1}(c),(d)), or \smallskip
\item  a pair of
non-parallel rays, directed by vectors with nonnegative
coordinates (see Figure \ref{nfig1}(e),(f)).
\end{itemize} (The labels $e_1$ and $e_2$ in Figure \ref{nfig1} denote the edges of $S({f})$ adjacent to
the point $(k,\ell)$, 
the symbol $\Delta$ designates the side of the depicted fragment
of the boundary on which the Newton polygon lies; we also note
that in cases (a), (c), (e), and (f), the rays drawn in bold may
merge to the same vertex, this does not affect our argument.)

\begin{figure}
\setlength{\unitlength}{0.6cm}
\begin{picture}(12,34)(0,0)
%
\thinlines \put(1,1){\vector(1,0){4}}\put(1,1){\vector(0,1){3}}
\put(1,6){\vector(1,0){4}}\put(1,6){\vector(0,1){4}}
\put(1,12){\vector(1,0){4}}\put(1,12){\vector(0,1){4}}
\put(1,18){\vector(1,0){4}}\put(1,18){\vector(0,1){4}}
\put(1,24){\vector(1,0){4}}\put(1,24){\vector(0,1){4}}
\put(1,30){\vector(1,0){4}}\put(1,30){\vector(0,1){4}}
%
\thicklines \put(1,3){\line(1,0){2}}\put(3,1){\line(0,1){2}}
\put(1,8){\line(1,0){1}}\put(2,8){\line(1,-1){2}}
\put(1,15){\line(1,-1){2}}\put(3,13){\line(1,-2){0.5}}
\put(1,21){\line(2,-1){2}}\put(3,20){\line(1,-1){1}}
\put(4,19){\line(1,-2){0.5}}\put(1,27){\line(1,-2){1}}
\put(2,25){\line(2,-1){2}}\put(1,33){\line(1,-2){0.5}}
\put(1.5,32){\line(1,-1){1}}\put(2.5,31){\line(2,-1){2}}
%
\put(0.7,2.8){$j$}
\put(0.7,8.8){$j$}\put(0.7,14.8){$j$}\put(0.7,20.8){$j$}
\put(0.7,26.8){$j$}\put(0.7,32.8){$j$}\put(2.9,0.5){$i$}\put(3.9,5.5){$i$}
\put(3.4,11.5){$i$}\put(4.4,17.5){$i$}\put(3.9,23.5){$i$}
\put(4.4,29.5){$i$}
%
\put(2.5,3.2){$(k,\ell)$}\put(2,8.2){$(k,\ell)$}\put(2,14.2){$(k,\ell)$}
\put(3.1,20.2){$(k,\ell)$}\put(1.9,24.1){$(k,\ell)$}\put(1.9,30.4){$(k,\ell)$}
\put(2.85,2.85){$\bullet$}
\put(1.85,7.85){$\bullet$}\put(1.85,13.9){$\bullet$}\put(2.85,19.9){$\bullet$}
\put(2.9,24.35){$\bullet$}\put(2.4,30.85){$\bullet$}
\put(2.85,12.9){$\bullet$}\put(3.85,18.9){$\bullet$}
\put(1.9,24.85){$\bullet$}\put(1.4,31.85){$\bullet$}
%
\put(3.1,19.2){$e_2$} \put(3.4,24.5){$e_2$}\put(3.5,30.6){$e_2$}
\put(1.4,7.6){$e_1$}\put(1.4,2.6){$e_1$}
\put(1.2,14.1){$e_1$}\put(1.8,20.1){$e_1$}\put(2.4,25){$e_1$}\put(2.1,31.6){$e_1$}
\put(2.5,2){$e_2$}\put(2.5,6.8){$e_2$}\put(2.2,13.1){$e_2$}
%
\put(1.5,1.5){$\Delta$} \put(1.5,6.5){$\Delta$}
\put(1.5,12.5){$\Delta$}\put(1.5,18.5){$\Delta$}\put(2.5,26){$\Delta$}
\put(3,32){$\Delta$}
%
\put(6,2){$\longleftrightarrow$}\put(6,7.5){$\longleftrightarrow$}
\put(6,13.5){$\longleftrightarrow$}\put(6,19.5){$\longleftrightarrow$}
\put(6,25.5){$\longleftrightarrow$}
\put(6,31.5){$\longleftrightarrow$}

%
\put(9,2){\line(0,1){2}}\put(10,1.5){\line(1,0){1.5}}
\put(9,7){\line(0,1){2}}\put(10,6.5){\line(1,1){2}}\put(9,13){\line(1,1){2}}
\put(10,12.5){\line(1,1){2}}\put(9,19){\line(1,2){1}}\put(10,18.5){\line(1,1){2}}
\put(8,24){\line(1,2){1}}\put(9,24){\line(1,2){1}}\put(8,30){\line(1,1){2}}\put(10,30){\line(1,2){1}}
\thinlines
\dashline{0.2}(8,2)(9,2)\dashline{0.2}(9,2)(10,1.5)\dashline{0.2}(10,1.5)(10,1)
\dashline{0.2}(8,7)(9,7)\dashline{0.2}(9,7)(10,6.5)\dashline{0.2}(10,6.5)(10,6)
\dashline{0.2}(8,13)(9,13)\dashline{0.2}(9,13)(10,12.5)\dashline{0.2}(10,12.5)(10,12)
\dashline{0.2}(8,19)(9,19)\dashline{0.2}(9,19)(10,18.5)\dashline{0.2}(10,18.5)(10,18)
\dashline{0.2}(8,27)(9,26)\dashline{0.2}(9,26)(10,26)\dashline{0.2}(10,26)(11,27)
\dashline{0.2}(9,33)(10,32)\dashline{0.2}(10,32)(11,32)\dashline{0.2}(11,32)(12,33)
%
\put(12,2){$U$}\put(12,7){$
U$}\put(12,13){$U$}\put(12,19){$U$}\put(11,25){$U$}\put(12,31){$U$}
%
\put(6.2,0){{\rm (f)}}\put(6.2,5){{\rm (e)}}\put(6.2,11){{\rm
(d)}}\put(6.2,17){{\rm (c)}}\put(6.2,23){{\rm
(b)}}\put(6.2,29){{\rm (a)}}
\end{picture}
\caption{Illustration to the proof of Theorem
\ref{nt1}}\label{nfig1}
\end{figure}
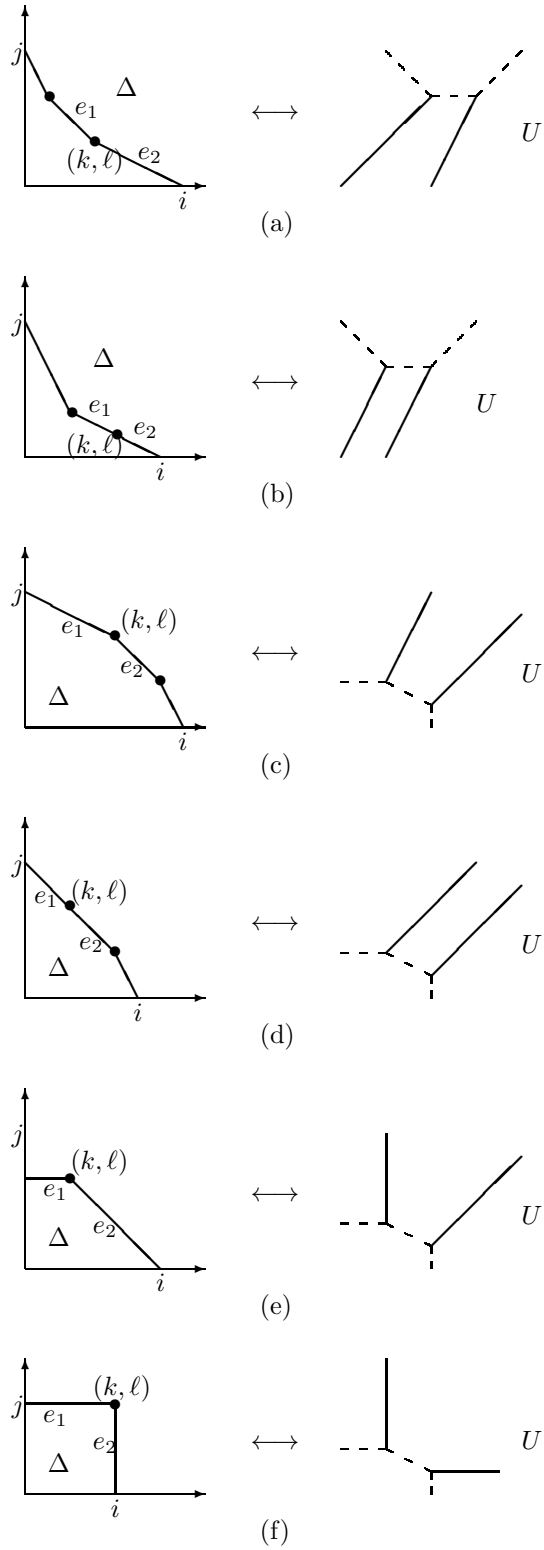

In all the situations described above, the tropical sums of points
lying on such pairs of rays sweep a two-dimensional domain in
$\Real^2$, contradicting the one-dimensionality of $P$.

Assume  $f$ does not contain a monomial $\a_i \TrP \var_1^i$ with
$i>0$. Since $f$ is not divisible by any monomial, $f$ should also
contain  a constant term $\a_0\in\Real$. This yields that $f$ has
no mixed monomials $\a_{k\ell}\TrP \var_1^k\TrP \var_2^\ell$, with
$k,\ell>0$, since otherwise, the Newton polygon $\Delta$ would
have a side with an  outer normal whose coordinates are nonzero
and having distinct signs -- a contradiction to the conclusion of
Step 2. Therefore, $f$ is simple.

{\it Step 4}. Suppose that $n\ge 3$, and $P=Z(f)$ with
$f\in\T[\Var]$.

Write $f$ as the sum of  essential monomials:
$f=\bigoplus_{\omega\in\Omega}M_\omega$, where $\Omega\subset\Z^n$
is finite. Assume that  $f$ has an essential monomial $M_\tau$,
$\tau\in\Omega$, depending on at least two variables, say
$\var_1,\var_2$. By definition, there are $c_1,\dots,c_n\in\Real$
for which
\begin{equation}M_\tau(c_1,\dots,c_n)>M_\omega(c_1,\dots,c_n)\quad\text{for
each}\quad\omega\in\Omega\backslash\{\tau\}\
.\label{ne11}\end{equation} Since a small variation of
$c_1,\dots,c_n$ does not violate (\ref{ne11}), we can take  these
numbers to be generic. For generality, the precise requirement is
as follows:  denoting by $\pr_{12}:\Z^n\to\Z^2$ the projection of
$\Z^n$ to the two first coordinates of $\Z^n$,  we rewrite the
polynomial $f$ as
$$f(\var_1,\dots,\var_n)=\bigoplus_{(k_1,k_2)\in\pr_{12}(\Omega)}\var_1^{k_1}\TrP
\var_2^{k_2}\TrP f_{k_1k_2}(\var_3,\dots,\var_n)\ ,$$
where
\begin{equation}f_{k_1k_2}(\var_3,\dots,\var_n)=\bigoplus_{(k_1,k_2,k_3,\dots,k_n)\in\Omega}
\a_{k_1 \cdots k_n}\TrP \var_3^{k_3}\TrP \cdots \TrP
\var_n^{k_n},\quad (k_1,k_2)\in\pr_{12}(\Omega)\
.\label{ne5}\end{equation}
Our demand is that, for each polynomial (\ref{ne5}), the values of
the monomials at $(c_3,\dots,c_n)$ be distinct. Geometrically,
this means that $(c_3,\dots,c_n)$ lies outside
$\bigcup_{(k_1,k_2)\in\pr_{12}(\Omega)}Z(f_{k_1k_2})\subset\Real^{n-2}$,
and that such a generic choice is always possible, since the
latter set is a finite polyhedral complex of dimension $n-3$ in
$\Real^{n-2}$.

Let $\Pi=\{\var_3=c_3, \dots ,\var_n=c_n\}\subset\Real^n$ be a plane
in $\Real^n$, and let  $g\in\T[\var_1,\var_2]$ be the polynomial
obtained from $f$ by substituting $c_3,\dots,c_n$ for
$\var_3,\dots,\var_n$, respectively. We claim that $P_2:=P\cap\Pi$
is the tropical set-curve in $\Pi$ given by $g$. Indeed, if
$(\u_1,\u_2)\in Z(g)$, then
\begin{equation}\u_1^{k_1}\TrP \u_2^{k_2}\TrP
f_{k_1k_2}(c_3,\dots,c_n)=\u_1^{\ell_1}\TrP \u_2^{\ell_2}\TrP
f_{\ell_1\ell_2}(c_3,\dots,c_n)\label{ne6}\end{equation}
for some $(k_1,k_2)\ne(\ell_1,\ell_2)\in\pr_{12}(\Omega)$, and
thus, $f$ has a pair of monomials reaching the value
$f(\u_1,u_2,c_3,\dots,c_n)$, namely $(\u_1,\u_2,c_3,\dots,c_n)\in
Z(f)$. On the other hand, if $(\u_1,\u_2,c_3,\dots,c_n)\in Z(f)$,
then the value $f(\u_1,\u_2,c_3,\dots,c_n)$ is attained by a pair
of monomials $M_{\omega'},M_{\omega''}$ of $f$, which in addition
must satisfy
$$\pr_{12}(\omega')=(k_1,k_2)\ne\pr_{12}(\omega'')=(\ell_1,\ell_2)$$ due to the choice of
$(c_3,\dots,c_n)$. Hence, equality (\ref{ne6}) is satisfied, which
means that $(\u_1,\u_2)\in Z(g)$.

The 2-plane $\Pi=\{\var_3=c_3,\dots,\var_n=c_n\}$ is a subgroup of
$(\Real^n,\TrS)$ isomorphic to $(\Real^2,\TrS)$, and therefore $P_2$
is an additive tropical set-curve in $\Real^2$.

To summarize, we have $\pr_{12}(\tau)=(k,\ell)$ with $k,\ell>0$
for monomial $M_\tau$ initially chosen to be essential. Then, the
monomial $N_{k\ell}=\var_1^k\TrP \var_2^\ell\TrP
f_{k\ell}(c_3,\dots,c_n)$ of $g$ is essential as well, since, due
to (\ref{ne11}), its value at $(c_1,c_2)$ is greater than that for
all the other monomials of $g$. As shown in Step 3, the sum $\hat
g$ of the essential monomials of $g$ must be a simple polynomial,
possibly multiplied by a monomial. Since $N_{k\ell}$ is essential
and depends both on $\var_1$ and $\var_2$, the polynomial $\hat g$
is divisible either by $\var_1$ or by $\var_2$. If $\hat g$ is
divisible by $\var_1$, then so is $g$. Indeed, otherwise, at least
one of the monomials of $g$ depends only on $\var_2$, and would
correspond to a vertex of the Newton polygon. Thus, it must be
essential (see Section \ref{nsec2}). Finally, notice that if $g$
is divisible by $\var_1$, then so is $f$, contradicting the
assumption of Step 1.

The proof of Theorem \ref{nt1} is completed.
\end{proof}

Another example of simple additive tropical sets is provided by
additive affine subspaces of $\Real^n$. Note that a hyperplane in
$\Real^n$ is a tropical set-hypersurface, since it can be defined by
a tropical binomial. Respectively, any rational affine subspace of
$\Real^n$ is a tropical set-variety defined by a number of tropical
binomials.

\begin{theorem}
An affine subspace $P\subset\Real^n$, parallel to a linear subspace
defined over $\Q$, is additive if and only if\; $P$ is simple.
\end{theorem}

\begin{proof}
As before, given an additive affine subspace $P\subset\Real^n$,
the task is to find simple tropical binomials that define $P$. In
view of Theorem \ref{nt1}, we may assume that $k=\dim (P) \le
n-2$. Choose a base $\bv_1,\dots,\bv_k$ of the linear space
parallel to $P$ and, without loss of generality, assume that the
first $k\times k$ minor of the coordinate matrix of
$\bv_1,\dots,\bv_k$ is nonsingular. Then $P$ projects onto a
hyperplane $P_i$ in the coordinate $(k+1)$-plane
$\Pi_i=\{\var_j=0,\ k<j\le n,\ j\ne i\}$, $i=k+1,\dots,n$. Using
Theorem \ref{nt1} again, we have $P_i=Z(f_i)\cap\Pi_i$, where
$f_i$ is a simple binomial for each $i=k+1,\dots,n$, and hence
$P=\bigcap_{i=k+1}^nZ(f_i)=Z(f_{k+1},\dots,f_n)$.
\end{proof}

\section{Additive tropical set-curves}
The treatment of additive tropical set-curves appears to be more
involved and delicate than one may expect. Our exposition appeals to
the natural idea of considering the projections of a given curve to
the coordinate planes and taking the intersection of the cylinders
built over all these projections. This intersection can be greater
than the original curve, and the central problem is then to remove
unnecessary pieces; this is what we are doing below. So, we proceed
as follows: first, we clarify several geometric properties of
additive tropical set-curves, later we construct some auxiliary
additive tropical sets, and, finally, we prove that additive
tropical set-curves are simple.

\subsection{Geometry of additive tropical set-curves}\label{secn1}
Let $U\subset\Real^n$ ($n\ge 2$) be an additive tropical set-curve.
Without loss of generality, we may assume that $U$ does not lie
entirely in any hyperplane $\var_j=\const$, $1\le j\le n$.

We denote the sets of the vertices and edges of $U$, respectively as
$U^0$ and $U^1$. Let us outline some useful geometric properties of
additive tropical set-curves. \medskip

\textbf{(i)} The directing vectors of the edges of $U$ may not have
a pair of coordinates having  distinct signs, since otherwise, along
the argument of Step 2 in the proof of Theorem \ref{nt1}, the sums
of points on such an edge would fill a two-dimensional domain.

We shall equip all the edges $e$ of $U$ with an orientation, taking
their (integral primitive) directing vectors $\ba(e)$ to have only
nonnegative coordinates. Note that this orientation agrees with the
order given by (\ref{ne3}). In addition, this orientation defines a
partial order in $U^1$ by letting $e\succ e'$ when $e$ and $e'$ have
a common vertex $\bu$, $e'$ coming to $\bu$, and $e$ emanating from
$\bu$. The poset $U^1$ has a unique maximal element, which is a ray
directed to $\Real^n_{\ge0}:=\{\v_1\ge0,\dots,\v_n\ge0\}$. Indeed,
otherwise, one would have two rays directed to $\Real^n_{\ge0}$, and
then, as was shown in Step 3 of the proof of Theorem \ref{nt1}, the
sums of points on such two edges would sweep a two-dimensional
domain.
\medskip

\textbf{(ii)} Let $\bu\in U^0$ and $C^1_{\bu}=\{e\in U^1\ :\ \bu\in
e\}$. As pointed above, in the notation of Section \ref{nsec3}, we
have
$$\ba_{\bu}(e)\in\Real^n_{\ge 0}\cup\Real^n_{\le 0}\quad\text{for all}\ e\in
U^1_{\bu},\quad\text{where}\
\Real^n_{\le0}=\{\v_1\le0,\dots,\v_n\le0\}\ .
$$ Furthermore, due to (\ref{ne2}), $U^1_{\bu}$ must contain at
least one edge $e$ with $\ba_{\bu}(e)\in\Real^n_{\ge 0}$ and at
least one edge $e'$ with $\ba_{\bu}(e')\in\Real^n_{\le 0}$. We also
claim  that $U^1_{\bu}$ contains precisely one edge $e$ with
$\ba_{\bu}(e)\in\Real^n_{\ge 0}$. Indeed, otherwise, the sums of
points on such two edges would sweep a two-dimensional domain. We
denote this edge by $e_{\bu}$. A similar reasoning shows that there
is at most one edge $e'$ with
\mbox{$\ba_{\bu}(e')\in\Real^n_{<0}:=\{\v_1<0,\dots,\v_n<0\}$}.
\medskip

\textbf{(iii)} Next, we notice that if
$\ba_{\bu}(e_{\bu})\in\{\v_i=0\}$, then $\ba_{\bu}(e')\in\{\v_i=0\}$
for all $e'\in U^1_{\bu}$. Indeed, if $\ba_{\bu}(e')$ has a nonzero
$i$-th coordinate for some $e'\in U^1_{\bu}$, then, due to
(\ref{ne2}), there should be some $\ba_{\bu}(e'')$, $e''\in
U^1_{\bu}$, with a positive $i$-th coordinate in contrary to
$\ba_{\bu}(e'')\in\Real^n_{\le 0}$ for all $e''\in
C^1_{\bu}\backslash\{e_{\bu}\}$. \medskip

\textbf{(iv)} Let $U_+$ denote the union of those edges $e\in U^1$
whose directing vectors satisfy
$$\ba(e)\in\Real^n_{>0}:=\{x_1>0,\dots,x_n>0\}\ .$$ We point out that
$U_+\ne\emptyset$, since it contains the maximal edge-ray $e\in
U^1$. Indeed, otherwise, by (iii) the whole tropical set-curve $U$
would lie in a hyperplane $x_i=\const$ -- contrary to the initial
assumption. Furthermore, due to (ii), $U_+$ must be connected and
homeomorphic either to $[0,\infty)$ or to $\Real$. We call $U_+$ the
\emph{spine} of the additive tropical set-curve $U$. \medskip

\textbf{(v)} Let $U^0_+=U_+\cap U^0:=\{\bu_1,\dots,\bu_m\}$ be the
set of vertices of $U$ that  lie on $U_+$. Pick
$i\in\{1,\dots,m\}$, and associate the set
$J(e)\subset\{1,\dots,n\}$ consisting of indices for which the
 coordinates of $\ba_{\bu_i}(e)$ are negative, to each edge $e\in U^1_{\bu_i}$ such that
$\ba_{\bu_i}(e)\in\Real^n_{\le 0}$. The additivity condition
implies that
\begin{itemize}\item the map $e\mapsto J(e)$,
restricted to $U^1_{\bu_i}$, is injective,
\item if $e_1,e_2\in
U^1_{\bu_i}$ emanate from $\bu_i$ in non-positive directions, then
either $J(e_1)\cap J(e_2)=\emptyset$, or there is an edge $e\in
U^1_{\bu_i}$ with $\ba_{\bu_i}(e)\in\Real^n_{\le 0}$ such that
$J(e)=J(e_1)\cap J(e_2)$.
\end{itemize}\medskip

\textbf{(vi)} Let  $U_i=\{\bu\in U\ :\ \bu\prec \bu_i\}$. This is
the part of the curve $U$ that  lies in the shifted orthant
$\bu_i+\Real^n_{\le 0}$. Since this orthant is a subsemigroup of
$\Real^n$, $U_i$ is an additive tropical set.

\subsection{Auxiliary additive tropical sets} Introducing  the
cone
$$ \Sig _0:=\Real^n_{\le0}\backslash\Real^n_{<0}=\{(u_1,\dots,u_n)\in\Real^n_{\le0}\ :\ u_1, \dots, u_n=0\}\ ,$$
and denoting  by $ \Sig _{\bu}$ the shift of $ \Sig _0$ to the
cone with vertex at $\bu\in\Real^n$, and by the results of Section
\ref{secn1}, we have
$$U\subset \wtU:=U_+\cup\bigcup_{\bu\in U^0_+} \Sig _{\bu}\ .$$
The cone $ \Sig _{\bu}$ divides $\Real^n$ into two components
which we denote by
$$\Inte( \Sig _{\bu})=\bu+\Real^n_{<0}\quad\text{and}\quad
\Ext( \Sig _{\bu})=\Real^n\backslash( \Sig _{\bu}\cup\Inte( \Sig
_{\bu}))\ .$$

The cone $ \Sig _0$ (and, respectively, each cone $ \Sig
_{\bu_i}$, $i=1,\dots,m$)  splits naturally into the disjoint
union of open cells, labeled by subsets
$J\subsetneq\{1,\dots,n\}$, and defined by
$$ \Sig _0(J)=\{(u_1,\dots,u_n)\in  \Sig _0\ :\ u_j<0\ \text{as}\ j\in J,\quad
u_j=0\ \text{as}\ j\not\in J\}\ .$$
Observing that
$$\overline{ \Sig _0(J)}=\bigcup_{K\subset J} \Sig _0(K)\ ,$$
we let
$$ \mcJ _i(U)=\{J\subsetneq\{1,\dots,n\}\ :\  \Sig _{\bu_i}(J)\cap
U\ne\emptyset\}\quad\text{and}\quad
 \Sig _{\bu_i}^U=\bigcup_{J\in{\mcJ}_i(U)}\overline{ \Sig _{\bu_i}(J)}\ ,$$ for each $i=1,\dots,m$,
and define
$$ \wtU_{\red}:=U_+\cup\bigcup_{\bu\in U_+^0} \Sig _{\bu}^U\ .$$
Note that $U\subset \wtU_{\red}\subset \wtU$ and that, for $n=2$,
$U= \wtU_{\red}$.

\begin{lemma}\label{lnew1}
$ \wtU$ and $ \wtU_{\red}$ are simple additive tropical sets.
\end{lemma}

\begin{proof}
We shall define $ \wtU$ and $ \wtU_{\red}$ by simple tropical
polynomials.

\medskip

\textbf{(1)} We first consider  $ \wtU$, and organize our argument
in a few steps.

\medskip

{\it Step 1}. Assume that $U_+$ is homeomorphic to $\Real$. We
intend to determine a (finite) set $\Phi$ consisting of simple
tropical polynomials in $n$ variables such that $\bigcap_{f\in
\Phi}Z(f)= \wtU$. In this step we show that
$\bigcap_{f\in\Phi}Z(f)\supset \wtU$.

Let
$$U^0_+=\{\bu_1,\dots,\bu_m\},\quad \bu_i=(u_{i1},\dots,u_{in}),\
i=1,\dots,m\ ,$$ with $u_{ij}<u_{kj}$ for all $1\le i<k\le m$,
$j=1,\dots,n$. The set $U_+$ contains $m+1$ edges, in order
$e_0\prec e_1=[\bu_1,\bu_2]\prec \cdots \prec
e_{m-1}=[\bu_{m-1},\bu_m]\prec e_m$, where $e_0$ and $e_m$ are
rays, whose primitive integral directing vectors are
$$\ba(e_i)=(a_{i1},\dots,a_{in})\in\Real^n_{>0},\quad i=0,\dots,m\
.$$
In particular, $u_{i+1,s}-u_{is}=a_{is}\mu_i$ for some $\mu_i>0$
for each $1\le i<m$ and $s=1, \dots ,n$.

Let $p_0,\dots,p_m$ and $b(i,j)$, $i=0, \dots ,m$, $j=1, \dots
,n$, be positive integers such that
\begin{enumerate}\item[(P1)] $p_i$ is divisible by $2a_{i1} \cdots a_{in}$,
$i=0,\dots,m$;
\smallskip
\item[(P2)] $b(0,j)\gg n$ and $b(i,j)-b(i-1,j)\gg n$
for all $i=1,\dots,m$, $j=1,\dots,n$, where
\mbox{$b(i,j):=p_i/a_{ij}$}.
\end{enumerate} By definition \begin{equation}a_{ik}b(i,k)=
a_{i\ell}b(i,\ell)\quad\text{for all}\quad k,\ell=1,\dots,n,\
i=1,\dots,m-1\ .\label{e13}\end{equation} Now we introduce the set
$\Phi\subset\T[\Var]$ of $n(n-1)/2$ simple tropical polynomials
$f_{k\ell}$, $1\le k<\ell\le n$, given by
\begin{equation}\label{flk}
\begin{array}{lll}
f_{k\ell} & =& \left(\bigoplus_{i=0}^m\left(\a^{k\ell}_{ik}\TrP
\var_k^{b(i,k)}\TrS \a^{k\ell}_{i\ell}\TrP
\var_\ell^{b(i,\ell)}\right)\right)
\\ [1mm] &&
\TrS\left(\bigoplus_{i=0}^{m-1}\bigoplus_{\renewcommand{\arraystretch}{0.6}
\begin{array}{c}
\scriptstyle{1\le j\le n}\\
\scriptstyle{j\ne k,\ell}
\end{array}}\left(\b^{k\ell}_{i+1,j}\TrP \var_j^{b(i+1,j)-j}\TrS \c^{k\ell}_{ij}\TrP
\var_j^{b(i,j)+j}\right)\right)\ ,
\end{array}
\end{equation}
 whose coefficients $\a^*_*$ are as specified below. The monomials of $f_{k\ell}$ correspond to the following integral
points:
\begin{enumerate}\eroman
\item $m+1$ points $P_{ki}=b(i,k)\eps_k$, $i=0, \dots ,m$, on the $k$-th
axis; \smallskip

\item $m+1$ points $P_{\ell i}=b(i,\ell)\eps_\ell$, $i=0, \dots ,m$, on
the $\ell$-th axis; \smallskip

\item $2m$ points $P^+_{ji}=(b(i,j)+j)\eps_j$, $i=0,
\dots ,m-1$, and $P^-_{ji}=(b(i,j)-j)\eps_j$, $i=1, \dots ,m$, on
the $j$-th axis for all $1\le j\le n$, $j\ne,k,\ell $ (here
$\eps_j$ denotes the $j$-th unit orthant).
\end{enumerate}

 The Newton polytope of $f_{k\ell }$ naturally splits
off the subpolytopes
\begin{equation}\Pi^i_{k\ell }=\conv\{P_{k,i-1},\ P_{ki}, \ P_{\ell ,i-1}, \ P_{\ell i},P^+_{j,i-1}, \ P^-_{ji},\
j\ne k,\ell \},\quad i=1, \dots ,m\ .\label{ne70}\end{equation}
Now we impose conditions on the coefficients of $f_{k\ell }$:
\begin{equation}
\begin{array}{lll}
\left(\a^{k\ell }_{i-1,k}\TrP
\var_k^{b(i-1,k)}\right)\Big|_{\bu_i} & = & \left(\a^{k\ell
}_{i-1,\ell }\TrP \var_\ell ^{b(i-1,\ell )}\right)\Big|_{\bu_i}
 =  \left(\a^{k\ell }_{ik}\TrP
\var_k^{b(i,k)}\right)\Big|_{\bu_i} \\[1mm]
  & =&  \left(\a^{k\ell }_{i\ell }\TrP
\var_\ell ^{b(i,\ell )}\right)\Big|_{\bu_i} \\[2mm]
& =& \left(\b^{k\ell }_{ij}\TrP
\var_j^{b(i,j)-j}\right)\Big|_{\bu_i}=\left(\c^{k\ell
}_{i-1,j}\TrP \var_j^{b(i-1,j)+j}\right)\Big|_{\bu_i}\
,\label{ne7}
 \end{array}
\end{equation} for all $i=1,\dots,m$, $1\le k<\ell \le n$,
$1\le j\le n$, $j\ne k,\ell $. We should check the consistency of
system (\ref{ne7}), since each of the coefficients $\a^{k\ell
}_{ik}$, $\a^{k\ell }_{i\ell }$, $i=1, \dots ,m-1$, enters two
equations in this system. The verification goes as follows:
(\ref{ne7}) reads as
\begin{equation}\a^{k\ell }_{i-1,k}+u_{ik}b(i-1,k)=\a^{k\ell }_{i-1,\ell }+u_{i\ell }b(i-1,\ell )=
\a^{k\ell }_{ik}+u_{ik}b(i,k)=\a^{k\ell }_{i\ell }+u_{i\ell
}b(i,\ell ) \label{nnne7}\end{equation} as $1\le i\le m$; then, we
have to show that
$$\a^{k\ell }_{ik}+u_{ik}b(i,k)=\a^{k\ell }_{i\ell }+u_{i\ell }b(i,\ell )\quad\Longrightarrow\quad
\a^{k\ell }_{ik}+u_{i+1,k}b(i,k)=\a^{k\ell }_{i\ell }+u_{i+1,\ell
}b(i,\ell )$$ as $1\le i<m$, or, equivalently, that
$$(u_{i+1,k}-u_{ik})b(i,k)=(u_{i+1,\ell }-u_{i\ell })b(i,\ell ),\quad
1\le i<m\ ,$$ which finally reduces to assumption (\ref{e13}). The
solutions of (\ref{ne7}) form a one-parametric family,  we pick
one of these solutions.

Now consider the truncation of $f_{k\ell }$ to one variable (i.e.,
the sum of monomials containing the only a chosen variable). From
(\ref{nnne7}) and property (P2) above, we derive
$$A^{k\ell }_{i+1,k}+ \al _{i+1,k}b(i+1,k)=A^{k\ell }_{ik}+ \al _{i+1,k}b(i,k)$$ $$\Longrightarrow
\quad A^{k\ell }_{i+1,k}+ \al _{ik}b(i+1,k)<A^{k\ell }_{ik}+ \al
_{ik}b(i,k)= A^{k\ell }_{i-1,k}+ \al _{ik}b(i-1,k)\ ,$$ which
immediately generalizes  to
\begin{equation}\begin{cases}\left(A^{k\ell }_{ik}\TrP
\var_k^{b(i,k)}\right)\Big|_{\bu_i}&>\left(A^{k\ell }_{sk}\TrP
\var_k^{b(s,k)}\right)\Big|_{\bu_i}\\ \left(A^{k\ell }_{i\ell
}\TrP \var_\ell ^{b(i,\ell )}\right)\Big|_{\bu_i}&>\left(A^{k\ell
}_{s\ell }\TrP \var_\ell ^{b(s,\ell
)}\right)\Big|_{\bu_i}\end{cases}\quad\text{when}\ |i-s|\ge 2\
.\label{enn1}\end{equation} Similarly, we have
\begin{equation}\begin{cases}\left(B^{k\ell }_{i+1,j}\TrP
\var_j^{b(i+1,j)-j}\right)\Big|_{\bu_i}&>\left(B^{k\ell
}_{s+1,j}\TrP \var_j^{b(s+1,j)-j}\right)\Big|_{\bu_i},\\
\left(C^{k\ell }_{ij}\TrP
\var_j^{b(i,j)+j}\right)\Big|_{\bu_i}&>\left(C^{k\ell }_{sj}\TrP
\var_j^{b(s,j)+j}\right)\Big|_{\bu_i},\end{cases}\quad\text{when}\
i\ne s\ .\label{enn2}\end{equation} Altogether, this means that
all the monomials of the considered truncation are essential.

The latter property and condition (\ref{ne7}) imply that the
subdivision $S(f_{k\ell })$ of the Newton polytope
$\Delta(f_{k\ell })$ contains the polytopes $\Pi^{k\ell }_i$,
$i=1,\dots,m$, defined by (\ref{ne70}). Furthermore, each polytope
$\Pi^{k\ell }_i$ is dual to the vertex $\bu_i$ (in the sense of
Section \ref{nsec2}), and the polytope's edges lying on the
coordinate axes are dual to the facets of the cone $ \Sig _{u_i}$.

Next, the edge $[P_{ki},P_{\ell i}]$ of the subdivision
$S(f_{k\ell })$ is dual to a convex $(n-1)$-dimensional polyhedron
in $Z(f_{k\ell })$ that contains either the point $\bu_1$ as
$i=0$, or the points $\bu_{i-1}$ and $\bu_i$ as $1\le i<m$, or the
point $\bu_m$ as $i=m$. In the case of $1\le i<m$, due to
convexity, the referred polyhedron contains the whole edge $e_i$
of $U_+$. In the case when $i=0$ or $m$, due to the orthogonality
of $[P_{ki},P_{\ell i}]$ to this polyhedron and to the edge $e_i$
(the latter orthogonality emerges from (\ref{e13})), the
hyperplane spanned by the polyhedron contains the edge $e_i$.
Moreover, the following comparison of monomials says that the
polyhedron itself contains $e_i$: due to (\ref{e13}) and
(\ref{ne7}), for $\bu=\bu_1-t\ba(e_0)\in e_0$, $t>0$, and any
$j\ne k,\ell $, $1\le j\le n$, one has
$$
\begin{array}{lll}
\left(C^{k\ell }_{0j}\TrP \var_j^{b(0,j)+j}\right)\Big|_{\bu} &= &
\left(C^{k\ell }_{0j}\TrP
\var_j^{b(0,j)+j}\right)\Big|_{\bu_1}-ta_{0j}b(0,j)-ta_{0j}j \\[1mm]
& = &\left(A^{k\ell }_{0k}\TrP
\var_j^{b(0,k)}\right)\Big|_{\bu_1}-ta_{0k}b(0,k)-ta_{0j}j
\\[1mm]
&=& \left(A^{k\ell }_{0k}\TrP
\var_j^{b(0,k)}\right)\Big|_{\bu}-ta_{0j}j
\\[1mm]
&< & \left(A^{k\ell }_{0k}\TrP \var_j^{b(0,k)}\right)\Big|_{\bu}\
,\end{array}$$ and similarly, for $\bu=\bu_m+t\ba(e_m)\in e_m$,
$t>0$, $$\left(B^{k\ell }_{mj}\TrP
\var_j^{b(m,j)-j}\right)\Big|_{\bu}<\left(A^{k\ell }_{mk}\TrP
\var_k^{b(m,k)}\right)\Big|_{\bu}\ .
$$
Thus, $\bigcap_{k,\ell }Z(f_{k\ell })\supset \wtU$.

\medskip

{\it Step 2}. Let us prove the inverse relation $\bigcap_{k,\ell
}Z(f_{k\ell })\subset \wtU$. More precisely, we have to show that
outside the cones $ \Sig _{\bu}$, $\bu\in U^0_+$, the ideal
generated by the polynomials $f_{k\ell }$, $1\le k<\ell \le n$,
defines a subset of $U_+$.

First, we introduce extra notation referring to the splitting of
each polynomial $f_{k\ell }$, $1\le k<\ell \le n$, into the
following (tropical) sum:
$$\begin{array}{lll}
f_{k\ell }&=&f^{(0)}_{k\ell }\oplus f^{(1)}_{k\ell }\oplus \dots
\oplus
f^{(m)}_{k\ell }\ ,\\[1mm]
f^{(0)}_{k\ell } & =& A^{0k}_{k\ell }\TrP \var_k^{b(0,k)}\oplus
A^{0\ell }_{k\ell }\TrP \var_\ell ^{b(0,\ell )}\oplus
\bigoplus_{j\ne k,\ell }C^{0j}_{k\ell }\TrP \var_j^{b(0,j)+j}\ ,\\[1mm]
f^{(m)}_{k\ell } & =& A^{mk}_{k\ell }\TrP \var_k^{b(m,k)}\oplus
A^{m\ell }_{k\ell }\TrP \var_\ell ^{b(m,\ell )}\oplus
\bigoplus_{j\ne
k,\ell }B^{mj}_{k\ell }\TrP \var_j^{b(m,j)-j}\ , \\[1mm]
f^{(i)}_{k\ell }&=& A^{ik}_{k\ell }\TrP \var_k^{b(i,k)}\oplus
A^{i\ell }_{k\ell }\TrP \var_\ell ^{b(i,\ell )}\oplus
\bigoplus_{j\ne k,\ell }\left(B^{ij}_{k\ell }\TrP
\var_j^{b(i,j)-j}\oplus C^{ij}_{k\ell }\TrP
\var_j^{b(i,j)+j}\right),\quad 1\le i<m\ .\end{array}$$

Let $\bu'=(u'_1, \dots ,u'_n)\in\bigcap_{k,\ell }Z(f_{k\ell
})\cap\Ext( \Sig _{\bu_m})$. Without loss of generality, we may
assume that $\bu_m=0$ and $f_{k\ell }(\bu_m)=0$, for all $1\le
k<\ell \le n$. Then, in particular,
$$f^{(m)}_{k\ell }=\var_k^{b(m,k)}\oplus
\var_\ell ^{b(m,\ell )}\oplus \bigoplus_{j\ne k,\ell }
\var_j^{b(m,j)-j}\ ,$$ and the point $\bu$ is such that $u'_i>0$
as $i$ belongs to a nonempty subset $J\subset\{1, \dots ,n\}$, and
$u'_j\le 0$ as $j\not\in J$. It follows  immediately from
(\ref{ne7}), (\ref{enn1}), and (\ref{enn2}) that, for any fixed
$1\le k<\ell \le n$,  the top degree monomials of $f_{k\ell }$ in
the variables $\var_i$, $i\in J$, take positive values at $\bu'$.
These values are greater than the values taken by the other
monomials in $\var_i$, $i\in J$,  at $\bu'$. Similarly,  the
monomials in $\var_j$, $j\not\in J$, take negative values at
$\bu'$. This means that the geometry of $Z(f_{k\ell })$ in $\Ext(
\Sig _{u_m})$ is determined by the top degree monomials of
$f_{k\ell }$,  i.e. by $f^{(m)}_{k\ell }$.

Next, we have 
$$\max_{r=1, \dots ,n}\var_r^{b(m,r)}\big|_{u'_r}=\var_i^{b(m,i)}\big|_{u'_i}>
\var_j^{b(m,j)}\big|_{u'_j},\quad i\in J',\ j\not\in J'\ ,$$ for
some set $J'\subset J$. Assuming that $J'\subsetneq\{1, \dots
,n\}$, we pick $k\in J'$ and $\ell \in\{1, \dots ,n\}\backslash
J'$, and obtain the following:
$$\var_k^{b(m,k)}\big|_{u'_k}=\var_i^{b(m,i)}\big|_{u'_i}>
\var_i^{b(m,i)-i}\big|_{u'_i}\quad\text{for all}\quad i\in
J'\backslash\{k\}\ ,$$
$$\var_k^{b(m,k)}\big|_{u'_k}>\var_i^{b(m,i)}\big|_{u'_i}\ \text{for all}\  i\in
J\backslash J',\quad\text{and}\ \var_k^{b(m,k)}\big|_{u'_k}>0\ge
\var_i^{b(m,i)-i}\big|_{u'_i}\ \text{for all}\ i\not\in J\ ;$$
this means that the value $f^{(m)}_{k\ell }(\bu')$ is attained by
a unique monomial,  i.e., $\bu'\not\in Z(f_{k\ell })$. Hence,
$J'=\{1, \dots ,n\}$, which, due to (\ref{e13}), implies that
$\bu'=\mu \ba(e_m)$ with $\mu>0$; that is $\bu'\in U_+$.

Let $\bu'\in\bigcap_{k,\ell }Z(f_{k\ell })\cap\Ext( \Sig
_{\bu_{m-1}})\cap\Inte( \Sig _{\bu_m})$, that is $\bu'=(u'_1,
\dots ,u'_n)$ such that $0>u'_i>u_{m-1,i}$, $i\in J$ and
$0>u_{m-1,j}\ge u'_j$, for $j\not\in J$, where $J\subset\{1, \dots
,n\}$ is some nonempty set. These relations, together with
equalities (\ref{ne7}) and inequalities (\ref{enn1}),
(\ref{enn2}), yield that
$$f^{(m-1)}_{k\ell }(\bu')>\max\{f^{(0)}_{k\ell }(\bu'), \dots ,f^{(m-2)}_{k\ell }(\bu'),f^{(m)}_{k\ell }(\bu')\}
\quad 1\le k<\ell \le n\ ;$$ nevertheless  the value
$f^{(m-1)}_{k\ell }(\bu')$ can be attained only by monomials which
depend on $\var_j$, $j\in J$.

In view of $\bu_m=0$, $f_{k\ell }(\bu_m)=0$, and equalities
(\ref{ne7}) for $i=m$, we then get
$$f^{(m-1)}_{k\ell }=
\var_k^{b(m-1,k)}\oplus \var_\ell ^{b(m-1,\ell )}\oplus
\bigoplus_{j\ne k,\ell }\left(B^{m-1,j}_{k\ell }\TrP
\var_j^{b(m-1,j)-j}\oplus \var_j^{b(m-1,j)+j}\right)\ .$$
Furthermore, due to equalities (\ref{ne7}) for $i=m-1$ and
inequalities (\ref{enn1}), (\ref{enn2}), we have
$$\var_i^{b(m-1,i)}\big|_{u'_i}>\var_i^{b(m-1,i)}\big|_{u_{m-1,i}}=f_{ij}(\bu_{m-1})\quad\text{for
all}\quad i\in J,\ j\ne i\ ,$$
$$\var_j^{b(m-1,j)}\big|_{u'_j}\le \var_j^{b(m-1,j)}\big|_{u_{m-1,j}}=f_{ij}(\bu_{m-1})\quad\text{for
all}\quad j\not\in J,\ i\ne j\ .$$ Hence,
\begin{equation}\max_{r=1, \dots ,n}\var_r^{b(m-1,r)}\big|_{u'_r}=\var_i^{b(m-1,i)}\big|_{u'_i}>
\var_j^{b(m-1,j)}\big|_{u'_j}\quad\text{as}\quad i\in J',\
j\not\in J'\label{enn3}\end{equation} for some nonempty set
$J'\subset J$. Suppose that $J'\subsetneq\{1, \dots ,n\}$, and
pick  $k\in J'$, $\ell \not\in J'$. Then
\begin{equation}\text{Eq. (\ref{enn3})}\quad\Longrightarrow\quad
\var_k^{b(m-1,k)}\big|_{u'_k}> \var_\ell ^{b(m-1,\ell
)}\big|_{u'_\ell }\ ,\label{enn4}\end{equation}
\begin{equation}\text{Eq. (\ref{enn3})}\ \&\ u'_i<0\quad
\Longrightarrow\quad \var_k^{b(m-1,k)}\big|_{u'_k}\ge
\var_i^{b(m-1,i)}\big|_{u'_i}>
\var_i^{b(m-1,i)+i}\big|_{u'_i},\quad i\ne k,\ell \
,\label{enn5}\end{equation}
$$\begin{array}{clll}
 \text{Eq. (\ref{ne7})}\ \&\ i\in J\backslash\{k\} & &
 \Longrightarrow  &
\var_k^{b(m-1,k)}\big|_{u'_k}-\left(B^{m-1,i}_{k\ell }\TrP
\var_i^{b(m-1,i)-i)}\right)\Big|_{u'_i} \\[1mm]
&& = & u'_kb(m-1,k)-(B^{m-1,i}_{k\ell }+u'_i(b(m-1,i)-i))\\[1mm]
&& = & (u'_k-u_{m-1,k})b(m-1,k)-(u'_i-u_{m-1,i})(b(m-1,i)-i)\\[1mm]
&& &  + \left(u_{m-1,k}b(m-1,k)-(B^{m-1,i}_{k\ell }+u_{m-1,i}(b(m-1,i)-i))\right)\\[1mm]
&& = & (u'_k-u_{m-1,k})b(m-1,k)-(u'_i-u_{m-1,i})(b(m-1,i)-i)\\[1mm]
&& &
+\left(\var_k^{b(m-1,k)}\big|_{u_{m-1,k}}-\left(B^{m-1,i}_{k\ell
}\TrP
\var_i^{b(m-1,i)-i}\right)\Big|_{u_{m-1,i}}\right)\\[1mm]
&&
\stackrel{(\ref{ne7})}{=} & (u'_k-u_{m-1,k})b(m-1,k)-(u'_i-u_{m-1,i})(b(m-1,i)-i)\\[1mm]
&& > & (u'_k-u_{m-1,k})b(m-1,k)-(u'_i-u_{m-1,i})b(m-1,i)\\[1mm]
&&
= & \left(\var_k^{b(m-1,k)}\big|_{u'_k}-\var_i^{b(m-1,i)}\big|_{u'_i}\right)\\[1mm]
&& &  +
\left(\var_k^{b(m-1,k)}\big|_{u_{m-1,k}}-\var_i^{b(m-1,i)}\big|_{u_{m-1,i}}\right)\
.
\end{array}$$

In view of $k\in J'$ and (\ref{enn3}), the former expression in
the last line is nonnegative. In its turn, the latter expression
vanishes;  this follows from (\ref{ne7}) and (\ref{nnne7}) (cf.
the verification of the consistency of system (\ref{ne7})
performed in Step 1). Hence,
\begin{equation}\var_k^{b(m-1,k)}\big|_{u'_k}>\left(B^{m-1,i}_{k\ell }\TrP
\var_i^{b(m-1,i)-i)}\right)\Big|_{u'_i}, \quad i\in
J\backslash\{k\}\ .\label{enn6}\end{equation} Finally, for
$i\not\in J$ and  $i\ne \ell $, one has
\begin{equation}
\begin{array}{lll}
\var_k^{b(m-1,k)}\big|_{u'_k} &
> & \var_k^{b(m-1,k)}\big|_{u_{m-1,k}} \\[1mm]
& =& f_{k\ell }(\bu_{m-1}) =  \left(B^{m-1,i}_{k\ell }\TrP
\var_i^{b(m-1,i)-i)}\right)\Big|_{u_{m-1,k}} \\[1mm]
&\ge &\left(B^{m-1,i}_{k\ell }\TrP
\var_i^{b(m-1,i)-i)}\right)\Big|_{u'_i}\
.\label{enn7}\end{array}\end{equation} Thus, the assumption
$J'\subsetneq\{1, \dots ,n\}$ together with
(\ref{enn4})-(\ref{enn7}) has led to the fact that the value
$f_{k\ell }(\bu')$ is attained by the unique monomial
$\var_k^{b(m-1,k)}$, namely $\bu'\not\in Z(f_{k\ell })$ -- a
contradiction. Hence, $J'=\{1, \dots ,n\}$, which, due to
(\ref{e13}), implies $\bu'-\bu_{m-1}=\lambda\ba(e_{m-1})$, that is
$\bu'\in e_{m-1}\subset U_+\subset \wtU$.

In the same manner we proceed further showing that if
$\bu'\in\bigcap_{k,\ell }Z(f_{k\ell })\cap\Ext( \Sig
_{\bu_{r-1}})\cap\Inte( \Sig _{\bu_r})$, than $\bu'\in
e_{r-1}\subset \wtU$, $r<m$, and then deducing the required
relation $\bigcap_{k,\ell }Z(f_{k\ell })\subset \wtU$.

\medskip

{\it Step 3}. In the remaining situation, when $U_+$ is
homeomorphic to $[0,\infty)$, we modify the preceding construction
in order to exclude any ray $e_0$ attached to the vertex $u_1$ and
directed to the negative infinity. Namely, in formula (\ref{flk})
for $f_{k\ell }$, we replace all the terms having exponents
$b(0,j)$, $j=1,\dots,n$, by a constant $A^{k\ell }_0$ which
satisfies condition (\ref{ne7}) for $i=0$. The equality
$\bigcap_{k,\ell }Z(f_{k\ell })= \wtU$ is then obtained in the
same way as in Steps 1 and 2 for the case $U_+\simeq\Real$.

\medskip

\textbf{(2)} To prove that $ \wtU_{\red}$ is simple, we extend the
ideal $I=\langle f_{k\ell }\ :\ 1\le k<\ell \le n\rangle$,
defining $ \wtU$, with the extra simple tropical polynomials
constructed below.

We assume that $U_+$ is homeomorphic to $\Real$. As in the
preceding situation,  in order cover the case of $U_+$
homeomorphic to $[0,\infty)$, the forthcoming construction should
be slightly modified; however, we skip this case.

Fix some  $i=1,\dots,m$, and choose a set
$K\subsetneq\{1,\dots,n\}$ such that $ \Sig _{u_i}(K)\not\subset
\wtU_{\red}$ (or, equivalently, $ \Sig _{\bu_i}(K)\cap
\wtU_{\red}=\emptyset$). Then, we shall construct a simple
polynomial $f_{i,K}$ such that the set $Z(f_{i,K})$ contains the
following:
\begin{itemize}
    \item the spine $U_+$, \smallskip
    \item all the cones $ \Sig _{\bu_k}$ for $1\le k\le m$,
$k\ne i$, \smallskip
    \item  and all the orthants $ \Sig _{\bu_i}(J)$ such that
$ \Sig _{\bu_i}(J)\subset \wtU_{\red}$, but $Z(f_{i,K})\cap \Sig
_{u_i}(K)=\emptyset$.
\end{itemize}
 Taking an appropriate $K\subsetneq\{1,
\dots ,n\}$ and  adding such simple polynomials for all $i=1,
\dots ,m$, we obtain the required ideal.

Further on, the required polynomial $f_{i,K}$ will be defined by
an explicit formula. Aiming to obtain uniform expressions, we
(formally) pick two extra vertices in $U_+$: a point $\bu_0\in
e_0\backslash\{\bu_1\}$ and a point $\bu_{m+1}\in
e_m\backslash\{\bu_m\}$. Accordingly, we add two more elements
$$b(m+1,j):=2b(m,j),\quad b(-1,j):=b(0,j)/2,\quad j=1,\dots,n\ ,$$
to the sequence $b(k,j)$, $0\le k\le m$, $1\le j\le n$, defined in
(P2).

To make clearer the construction and properties of $f_{i,K}$, we
start with an auxiliary polynomial which can be viewed as a
simplified version of the polynomial $f_{k\ell }$, as introduced
in the preceding step,
\begin{equation}f=\sum_{\renewcommand{\arraystretch}{0.6}
\begin{array}{c}
\scriptstyle{-1\le k\le m+1}\\
\scriptstyle{1\le j\le n}
\end{array}}A_{k,j}\TrP \var_j^{b(k,j)}\ ,\label{nne8}\end{equation} where, for all
$k=0, \dots ,m+1$, $$\left(A_{k-1,j}\TrP
\var_j^{(k-1,j)}\right)\Big|_{\bu_k}=\left(A_{k,l}\TrP
\var_j^{(k,\ell )}\right)\Big|_{\bu_k},\quad j,\ell =1, \dots ,n\
.$$

The argument of Step 1  implies immediately that: the conditions
imposed on the coefficients $A_{kj}$ are consistent; all the
monomials of $f$ are essential; and the subdivision $S(f)$ of
$\Delta(f)$ consists of the polytopes (cf. (\ref{ne70}))
\begin{equation}\Pi^k=\conv\{P_{jk},\ P_{\ell ,k-1},\
j,\ell =1, \dots ,n\},\quad k=0, \dots ,m+1\
,\label{nne11}\end{equation} which are dual to the vertices
$\bu_0, \dots ,\bu_{m+1}$ (here $P_{kj}=b(k,j)\eps_k$ are the
vertices of the polytopes (\ref{ne70})). In addition, one obtains
that $Z(f)\supset \wtU$.

Next we modify formula (\ref{nne8}). Pick an element $j_0\in\{1,
\dots ,n\}\backslash K$ and define the desired polynomial
$f_{i,K}$ to be
$$f_{i,K}=\sum_{\renewcommand{\arraystretch}{0.6}
\begin{array}{c}
\scriptstyle{-1\le k\le m+1}\\
\scriptstyle{1\le j\le n}
\end{array}}\hat A_{k,j}\TrP \var_j^{\hat b(k,j)}\ ,$$ where $\hat
b(k,j)=b(k,j)$ for all $k=-1, \dots ,m+1$, $j=1, \dots ,n$, except
for the cases $$\hat b(i,j)=b(i,j)+1,\ j\not\in K,\quad \hat
b(i-1,j)=b(i-1,j)-1,\ j\not\in K\cup\{j_0\}\ ,$$ while, for all
$k=0, \dots ,m+1$, $k\ne i$, the coefficients $\hat A_{kj}$
satisfy the condition
\begin{equation}\left(A_{k-1,j}\TrP
\var_j^{(k-1,j)}\right)\Big|_{\bu_k}=\left(A_{k,l}\TrP
\var_j^{(k,\ell )}\right)\Big|_{\bu_k},\quad j,\ell =1, \dots ,n\
,\label{nne9}\end{equation} and the new condition
\begin{equation}\left(A_{i-1,j}\TrP
\var_j^{(i-1,j)}\right)\Big|_{\bu_i}=\left(A_{i,\ell }\TrP
\var_j^{(i,\ell )}\right)\Big|_{\bu_i},\quad j\in K\cup\{j_0\},\
\ell \in K\ .\label{nne10}\end{equation}

Again, as in  Step 1, from (\ref{e13}), we derive the consistence
of conditions (\ref{nne9}) and (\ref{nne10}), as well as the
inequalities
$$\left(\hat A_{is}\TrP \var_s^{b(i,s)+1}\right)\Big|_{\bu_i}<\left(\hat A_{ij}\TrP
\var_j^{b(i,j)}\right)\Big|_{\bu_i},\quad j\in K,\ s\not\in K\ ,$$
$$\left(\hat A_{i-1,\ell }\TrP \var_\ell ^{b(i-1,\ell )-1}\right)\Big|_{\bu_i}<\left(\hat A_{i-1,s}\TrP
\var_s^{b(i-1,s)}\right)\Big|_{\bu_i},\quad s\in K\cup\{j_0\},\
\ell \not\in K\cup\{j_0\}\ .$$ These inequalities yield that all
the monomials of $f_{i,K}$ are essential, and that the subdivision
$S(f_{i,K})$ of the Newton polytope $\Del(f_{i,K})$ contains the
$n$-dimensional polytopes
$$\hat\Pi^k=\conv\{\hat P_{jk},\ \hat P_{\ell ,k-1},\ j,\ell =1, \dots
,n\},\quad k=0, \dots ,m+1,\ k\ne i,\quad \hat
P_{jk}=b(k,j)\eps_j\ ,$$ ($P_{jk}=b(k,j)\eps_j$ being the vertices
of the polytopes (\ref{ne70})), and the polytope
$$\hat\Pi_i=\conv\{\hat P_{j,i-1},\ \hat P_{\ell i},\ j\in
K\cup\{j_0\},\ \ell \in K\ .$$ Further on, the above polytopes
$\Pi^k$, $k=0, \dots i-1,i+1, \dots ,m+1$, are dual to the
vertices $\bu_0, \dots ,\bu_{i-1},\bu_{i+1}, \dots ,\bu_{m+1}$ of
$U_+$, and the polytope $\Pi_i$ is dual to a face of $Z(f_{i,K})$
which passes through $\bu_i$.

So, we immediately decide that $Z(f_{i,K})$ contains: the part of
$U_+$, preceding the vertex $\bu_{i-1}$, the part of $U_+$,
following the vertex $\bu_{i+1}$, and  all the cones $ \Sig
_{\bu_k}$, $k=0, \dots ,m+1$, $k\ne i$. Now, we observe that $K$
contains at least two elements. Indeed, otherwise, if $K=\{j_1\}$,
then the vectors $\ba_{\bu_i}(e)$, with $e\in U^1_{u_i}$, oriented
to $\Real^n_{\le0}$, would lie in the same hyperplane
$\{\var_{j_1}=u_{ij_1}\}$, and thus could not be balanced by a
vector $\ba_{\bu_i}(e_{\bu_i})\in\Real^n_{>0}$, which contradicts
(\ref{ne2}). So, if $j_1,j_2\in K$, then the $(n-1)$-face of
$Z(f_{i,K})$, dual to the edge $[P_{j_1,i-1},P_{j_2,i-1}]$,
contains the points $\bu_{i-1}$ and $\bu_i$, and hence, also
contains the edge $e_i$ of $U_+$. Similarly, the $(n-1)$-face of
$Z(f_{i,K})$, dual to the edge $[P_{j_1i},P_{j_2i}]$, contains the
points $\bu_i$ and $\bu_{i+1}$, and hence, also contains the edge
$e_{i+1}$ of $U_+$. That is $U_+\subset Z(f_{i,K})$.

Next we verify that $ \Sig _{\bu_i}(J)\subset \wtU_{\red}$ implies
$ \Sig _{\bu_i}(J)\subset Z(f_{i,K})$. Indeed, if $ \Sig
_{u_i}(J)\subset  \Sig _{\bu_i}^U$, then by construction
$J\not\supset K$. Hence, there exists $s\in K\backslash J$, and
thus the value $f_{i,K}(\bu_i)$ is attained (among others) by the
two monomials $\hat A_{is}\TrP \var_s^{b(i,s)}$ and $\hat
A_{i-1,s}\TrP \var_s^{b(i-1,s)}$. Then, the $(n-1)$-dimensional
orthant
$$\{\var_s=u_{is},\quad \var_j\le u_{ij},\ 1\le j\le n,\ j\ne
s\},$$ having the vertex $\bu_i$, is contained in $Z(f_{i,K})$,
and in its turn contains $ \Sig _{\bu_i}(J)$.

The last task is to check that $ \Sig _{\bu_i}(K)\cap
Z(f_{i,K})=\emptyset$. To show this, we note that the polynomial
$f_{i,K}$ is constant along
$$ \Sig _{\bu_i}(K)=\{\var_j= \al _{ij},\ j\not\in K,\quad
\var_l<u_{i\ell},\ \ell \in K\}$$ and its value is attained  only
by the monomial $\hat A_{ij_0}\TrP \var_{j_0}^{b(i,j_0)}$.

This completes the proof of Lemma \ref{lnew1}.
\end{proof}

The construction of the ideals defining $ \wtU$ and $\wtU_{\red}$
depends on the choice of the parameters $p_0,\dots,p_m$. Next we
define these parameters so that Proposition \ref{p1} below holds
true.

Given two strictly increasing sequences
$\overline\xi=\{\xi_1,\dots,\xi_r\}$ and
$\overline\eta=\{\eta_1,\dots,\eta_r\}$ of real numbers, we say
that $\overline\eta$ is \emph{$\overline\xi$-convex} if
\begin{equation}\frac{\eta_k-\eta_{k-1}}{\xi_k-\xi_{k-1}}<\frac{\eta_{k+1}-\eta_k}{\xi_{k+1}-\xi_k}\quad
\text{for all}\quad k=2,\dots,r-1\ .\label{e14}\end{equation}

\begin{proposition}\label{p1}
In the above notation, let
$\overline\xi^{(k)}=\{\xi^{(k)}_1,\dots,\xi^{(k)}_m\}$,
$k=1,\dots,s$, be an arbitrary strictly increasing sequences of
real numbers. Then, there  exist integers $p_0,\dots,p_m$,
satisfying conditions (P1), (P2) from the first part of the proof
of Lemma \ref{lnew1}, such that, for each generator $f$ of the
defining ideal of $ \wtU_{\red}$ and for every $k=1,\dots,s$, the
(strictly increasing) sequence $f(\bu_1),\dots,f(\bu_m)$ is
$\overline\xi^{(k)}$-convex.
\end{proposition}

We leave the proof of this elementary statement to the reader,
remarking only that one should choose the sequence $p_0,\dots,p_m$
which grows sufficiently quickly.

\subsection{Remark on plane additive tropical set-curves}\label{secn4} The above geometric
treatment, as well as the algebraic one, becomes quite transparent
in the case of  additive tropical plane set-curves.

Geometrically, one obtains an additive tropical set-curve
$U\subset\Real^2$ from its spine $U_+$ by attaching to each vertex
$\bu_i$, $1\le i\le m$, one or two negatively directed horizontal
and vertical rays. Furthermore, if $\bu_1$ is the minimal point of
the spine $U_+$ (i.e., $U_+\simeq[0,\infty)$), then we call $\bu_1$
a \emph{terminal vertex} of $U$. In particular, if $\bu_1$ is
terminal, then it is a common vertex of a horizontal and a vertical
negatively directed rays of $U$.

By Theorem \ref{nt1}, such a set-curve $U$ can be defined by one
simple tropical polynomial. Furthermore, Lemma \ref{lnew1} provides
a family of such polynomials with parameters $p_0,\dots,p_m$
subjecting  to conditions (P1), (P2). Keeping the property declared
in Proposition \ref{p1}, we claim that one can vary these parameters
and gets the following additional property:

\begin{proposition}\label{p2}
In the above notation, assume that the point $\bu_m\in U_+^0$ is a
common vertex of a horizontal and a vertical negatively directed
rays. Then, for any polynomial $f(\var_1,\var_2)$ constructed for
$U$ as in the proof of Lemma \ref{lnew1}, keeping the values
$f(\bu_1),\dots,f(\bu_m)$ and the set $Z(f)$ unchanged, one can
make
\begin{equation}p_m\gg p_k\quad\text{and}\quad p_0\ll p_k\quad \text{for all}\quad
k=1,\dots,m-1\ .\label{e15}\end{equation}
\end{proposition}

Again, the proof is an easy exercise left to the reader. We only
observe, that, if $\bu_1$ is not terminal then $p_0=0$ satisfies
the requirements of the proposition. Also, in such a variation of
$p_0$ and $p_m$, the convexity property required in Proposition
\ref{p1} persists, since it depends only on $f(u_1),\dots,f(u_m)$.

\subsection{Simplicity of spacial additive tropical set-curves}

\begin{theorem}\label{thm:end}
A tropical set-curve $U\subset\Real^n$, where $n=2$ or $3$, is
additive if and only if it is simple.
\end{theorem}

\begin{proof} In view of Corollary \ref{thm:powerOfPoly} and
Theorem \ref{nt1}, it remains to prove that an additive tropical
set-curve $U\subset\Real^3$ is simple.

Pick a number $i=1,\dots,m$ and an element $J\in{\mcJ}_i(U)$, and
consider the set $U_{i,J}=U\cap\overline{ \Sig _{u_i}(J)}$. We
intend to construct a pair of simple polynomials, denoted $F$ and
$F'$,  for which
\begin{equation}Z(F)\cap Z(F')\supset U\quad\text{and}\quad Z(F)\cap Z(F')
\cap\overline{ \Sig _{u_i}(J)}=U_{i,J}\
.\label{nne73}\end{equation} Then, ranging over all $i=1,\dots,m$
and over all $J\in \mcJ _i(U)$, and adding all the newly acquired
polynomials to the ideal of $ \wtU_{\red}$, we obtain the desired
simple ideal defining $U$.

The case of $\#J=1$ is easy. Indeed, $U_{i,J}$ is just the ray
parallel to one of the coordinate axes, say $\var_1$-axis,
emanating from $\bu_i$, and pointing to $-\infty$. We project $U$
to the $(\var_1,\var_2)$-plane and obtain an additive tropical
plane curve $V$, which is defined by a simple polynomial
$F(\var_1,\var_2)$. It is then clear that $Z(F)\cap\overline{ \Sig
_{u_i}(J)}=U_{i,J}$.

So, it remains to consider the case when  $\# J=2$; thus, from now
on we assume $J=\{1,2\}$. We identify $\Real^2$ with the plane
$\{\var_3=0\}\subset\Real^3$ and introduce the natural projection
$\pi:\Real^3\to\Real^2$. Our strategy is to construct the
polynomials $F$ and $F'$ to be of the form $f(\var_1,\var_2)\TrS
g(\var_3)$, where $f(\var_1,\var_2)$ is a simple polynomial
defining a certain modification of the projection $\pi(U)$ in
$\Real^2$, and $g(\var_3)$ provides a correction of the
intersections of $Z(f)$ with the quadrants $ \Sig _{u_k}(J)$,
$k=1, \dots ,m$.

We proceed further in several steps.

\medskip

{\it Step 1}. In general, $\pi(\overline{ \Sig
_{\bu_i}(J)})\cap\pi(U)$ is greater than $\pi(\overline{ \Sig
_{\bu_i}(J)})\cap\pi(U_{i,J})$. In this step, we shall decide
which  parts of $U$ contribute to $\pi(\overline{ \Sig
_{\bu_i}(J)})\cap\pi(U)$ beyond $\pi(U_{i,J})$, and which parts do
not.

The set $V_{i,J}=\pi(U_{i,J})$ is an additive  tropical plane set,
which can be extended up to an additive tropical set-curve $\hat
V_{i,J}$ by attaching a ray with vertex $\pi(\bu_i)$ directed to
$\Real^2_{>0}$. If $\hat V_{i,J}$ has a terminal vertex
$\bv=(v_1,v_2)$ (which then must differ from $\bu_i$, since $ \Sig
_{\bu_i}(J)\cap U\ne \emptyset$), we let
$$Q_{i,J}=\{(x_1,x_2)\in\Real^2\ : \max_{j=1,2}(x_j-u_{ij})\le0\le\max_{j=1,2}(x_j-v_j)\}\
,$$ and otherwise we set $$Q_{i,J}=\{(x_1,x_2)\in\Real^2\ :\
\max_{j=1,2}(x_j-u_{ij})\le 0\}\ .$$ Geometrically,  in the latter
case, $Q_{i,J}$ is just a shifted negative quadrant, while in the
former case, $Q_{i,J}$ is the closed difference of two such
quadrants, one lies inside the interior of the other. Moreover,
$V_{i,J}\subset Q_{i,J}$, and $Q_{i,J}$ is the minimal figure of the
given shape, containing $V_{i,J}$.

We claim that: \begin{itemize}

\item for each $k>i$ and
$K\subset\{1,2,3\}$ such that $K\not\supset J$, one has
$\pi(U_{k,K})\cap Q_{i,J}=\emptyset$, \smallskip

\item $\pi(U_i)\cap
Q_{i,J}=V_{i,J}$, where $U_i=\{\bu\in U\ :\ \bu\prec \bu_i\}$.
\end{itemize}
The first relation is easy: if $l\in J\backslash K$, then any
point of $U_{k,K}$ has the $l$-th coordinate $u_{kl}>u_{il}$, and
hence its $\pi$-projection lies outside $Q_{i,J}$. To prove the
second relation, we note that for any point $\bu=(u_1,u_2,u_3)\in
U_i$ for which  $\pi(\bu)\in Q_{i,J}$, one has $u_3\le u_{i3}$,
and there always exists a point $\bv=(u'_1,u'_2,u_{i3})\in
V_{i,J}$ such that $u'_1\le u_1$ and $u'_2\le u_2$ (the latter
property is evident if $\hat V_{i,J}$ has no terminal vertex, and
it follows from (vii) in Section \ref{secn1} if $\hat V_{i,J}$ has
a terminal vertex). Then,
$$\bu\TrS \bv=(u_1,u_2,u_{i3})\in V_{i,J}\ ,$$ which, in particular, yields
$\pi(\bu)=\pi(\bu\TrS \bv)$.

Thus, the only parts of $U$, whose $\pi$-projections may
contribute to $\pi(U)\cap Q_{i,J}$ beyond $\pi(U_{i,J})$, are
$U_{k,J}$, with $k>i$.

\medskip

{\it Step 2}. Assume that $\hat V_{i,J}$ has a terminal vertex
$\bv_1$ as appears in Step 1. In this situation we shall construct
just one required polynomial $F=F'$;  we start by constructing the
part $f(\var_1,\var_2)$ of $F$.

Pick  a point $\bu_{m+1}=(u_{m+1,1},u_{m+1,2},u_{m+1,3})$ on the
ray $e_m\subset U$, and attach to it the three rays parallel to
the coordinate axes and pointing to $-\infty$. The newly obtained
set $\hat U$ is again a tropical additive curve. Now, let
$$W_{i,J}=\hat U\cup\bigcup_{k>i}\overline{ \Sig _{\bu_k}(J)}
\backslash\bigcup_{k>i} \Sig _{\bu_k}(J)\ .$$ Geometrically,
$W_{i,J}$ is obtained from $\hat U$ as follows: for each vertex
$\bu_k$, with $k>i$, we delete the part of $U$ attached to $\bu_k$
and contained in the quadrant $Q_{k,J}$ and, instead, we add two
negatively directed rays emanating from $\bu_k$ and parallel to
the $\var_1$-axis and  to the $\var_2$-axis, respectively. It is
easy to see that $W_{i,J}$ is an additive tropical curve. Hence,
$\pi(W_{i,J})\subset\Real^2$ is an additive
 tropical plane curve. The results of Step 1 imply that
\begin{equation}\pi(W_{i,J})\cap Q_{i,J}=U_{i,J}\
.\label{nne74}\end{equation}

Note that the points $\pi(\bv_1)$ and
$\pi(\bu_1),\dots,\pi(\bu_m), \pi(\bu_{m+1})$ belong to the spine
of $\pi(W_{i,J})$. Let them be ordered (cf. \eqref{ne3}) as
follows
$$\pi(\bu_1)\prec \cdots \prec\pi(\bu_s)\prec\pi(\bv_1)\prec\pi(\bu_{i+1})\prec
\cdots \prec\pi(\bu_{m+1})\ ,$$ where $0\le s<i$ and $s$ is the
maximal possible index satisfying this ordering with
$\pi(\bu_s)\ne\pi(\bv_1)$. Due to Propositions \ref{p1} and
\ref{p2}, we may assume that $\pi(W_{i,J})$ is defined by a simple
polynomial $f(\var_1,\var_2)$ satisfying the following condition:
\\ the sequence
$$f(\bu_1)< \cdots <f(\bu_s)<f(\bv_1)<f(\bu_{i+1})< \cdots<f(\bu_{m+1})$$
is convex with respect to the sequence $u_{13}< \cdots
<u_{s3}<u_{i3}< \cdots <u_{m3}<u_{m+1,3}$, and relation
(\ref{e15}) holds true as well.

\medskip

{\it Step 3}. Now,  we define the polynomial
\begin{equation}F(\var_1,\var_2,\var_3)=f(\var_1,\var_2)\TrS g(\var_3),
\quad \text{with }  g(\var_3)=\bigoplus_{k=0}^{m+1}A_k\TrP
\var_3^{c_k}\ ,\label{e8}\end{equation}
whose parameters $A_k$ and $c_k$, $k=0,\dots,m+1$, satisfy the
following conditions:
\begin{enumerate}\item[(a)] $c_{m+1}=c_m+1$; \smallskip

\item[(b)] for
$i<k\le m$,
$$c_k=\frac{f(\bu_{k+1})-f(\bu_k)}{u_{k+1,3}-u_{k3}},\quad
\left(A_k\TrP \var_3^{c_k}\right)\Big|_{\bu_{k+1}}=
\left(A_{k+1}\TrP
\var_3^{c_{k+1}}\right)\Big|_{\bu_{k+1}}=f(\bu_{k+1})\ ;$$

\item[(c)] for the $i$-th monomial,
$$c_i=\frac{f(\bu_{i+1})-f(\bv_1)}{u_{i+1,3}-u_{i3}},\quad
\left(A_i\TrP \var_3^{c_i}\right)\Big|_{\bu_{i+1}}=
\left(A_{i+1}\TrP
\var_3^{c_{i+1}}\right)\Big|_{\bu_{i+1}}=f(\bu_{i+1})\ ,$$
$$\left(A_i\TrP \var_3^{c_i}\right)\Big|_{\bu_i}=f(\bv_1)\ ;$$

\item[(d)] $c_{i-1}=c_i-1$, and $\left(A_{i-1}\TrP
\var_3^{c_{i-1}}\right)_{\var_3=u_{i3}-\eps}=\left(A_i\TrP
\var_3^{c_i}\right)_{\var_3=u_{i3}-\eps}$, where $\eps>0$ is small
(we specify this later);\smallskip

 \item[(e)] for $s<k<i-1$,
$$c_k=c_{k+1}-1,\quad\left(A_k\TrP
\var_3^{c_k}\right)\Big|_{\bu_{k+1}}=\left(A_{k+1}\TrP
\var_3^{c_{k+1}}\right)\Big|_{\bu_{k+1}}\ ;$$

 \item[(f)]
$\left(A_s\TrP
\var_3^{c_s}\right)\Big|_{\bu_{s+1}}=\left(A_{s+1}\TrP
\var_3^{c_{s+1}}\right)\Big|_{\bu_{s+1}}$; \smallskip

\item[(g)] for
$0\le k<s$,
$$\left(A_k\TrP \var_3^{c_k}\right)\Big|_{\bu_{k+1}}=\left(A_{k+1}\TrP
\var_3^{c_{k+1}}\right)\Big|_{\bu_{k+1}}=f(\bu_{k+1})\ ;$$
\item[(h)] $c_0=c_1-1$.
\end{enumerate} We observe that these relations uniquely determine
the values of the arguments  $A_k$ and $c_k$, $k=0,\dots,m+1$, out
of the values $f(\bu_1),\dots,f(\bu_{m+1}),f(\bv_1)$. Multiplying
the parameters $p_0,\dots,p_m$ in the construction of the
polynomial $f$ by a suitable natural number, we multiply the
values of $f$ by that number, and thus we can achieve the
integrality of the exponents $c_k$ in the above definition.

Due to the assumed convexity property of the values of $f$,  each
monomials of $g$ is essential (here we specify the value of
$\eps$, taking into account that for $\eps=0$ all the monomials of
$g$ appear to be essential).

Let us verify that $Z(F)\supset U$. Observing that
$$Z(F)=\{f=g\} \ \cup \ (Z(f) \cap  \{f\ge
g\}) \ \cup \ (Z(g)\cap\{f\le g\})\ ,$$
we note that, for any point $\bu \in\{f=g\}$, the set $\{f\ge g\}$
contains the negative ray with vertex $\bu$, parallel to the
$\var_3$-axis, and that the set $\{f\le g\}$ contains the negative
quadrant with vertex $\bu$, parallel to the
$(\var_1,\var_2)$-plane. Notice also that $W_{i,J}\subset Z(f)$.
Then:

\begin{enumerate}
    \item  For any $k>i$, due to relations (b), (c), and the
construction in the proof of Lemma \ref{lnew1}, the value
$F(\bu_{k})$ is attained by the pair of monomials
$A_{k-1}\TrP \var_3^{c_{k-1}}$ and $A_k\TrP \var_3^{c_k}$ of
$g(\var_3)$; the same value $F(\bu_k)$ is also attained by some
four monomials of $f(\var_1,\var_2)$, since by construction the
plane tropical curve $\pi(W_{i,J})$ has four edges incident to its
vertex $\pi(\bu_k)$, two of them with positive slopes, and the two
others being negatively directed vertical and horizontal rays. It
is then easy to derive that $Z(F)\supset  \Sig _{\bu_k}$.
\smallskip

\item Since $f(\bu_{m+1})=g(\bu_{m+1})$ and $p_{m+1}\gg\max_lf(\bu_l)$
(cf. Proposition \ref{p2} and relation (a) above), and the
polynomial $g$ is linear in the half-space $\{\var_3\ge \al
_{m+1,3}\}$, we decide that $f(\bu)\ge g(\bu)$ along the ray
$e_{\bu_{m+1}}$, and hence $e_{\bu_{m+1}}\subset Z(F)$. Similarly,
$f(\bu_1)=g(\bu_1)$ and $p_0\ll\min_lf(\bu_l)$ (cf. Proposition
\ref{p2} and relation (h) above), and hence the ray of $\hat U$,
emanating from $\bu_1$ and proceeding  to $\Real^n_{<0}$ is
contained in $Z(F)$. \smallskip

\item By construction, for
$k>i$, the values $f(\bu_k)$ and $f(\bu_{k+1})$  are attained by
the same monomial, and the same also holds for $g$. Hence, due to
the linearity of $f$ and $g$ along the segment
$[\bu_k,\bu_{k+1}]$, this segment is contained in $Z(F)$. In the
same way, when $1\le k<s$ we have  $[\bu_k,\bu_{k+1}]\subset
Z(F)$.
\smallskip

\item Since $g(\bu_{i+1})=f(\bu_{i+1})$ and
$g(\bu_i)=f(\bv_1)<f(\bu_i)$, we derive that the segment
$[\bu_i,\bu_{i+1}]$ lies in the domain $\{f\ge g\}$, and hereby is
contained in $Z(F)$. \smallskip

\item We have $g(\bu_s)=f(\bu_s)$ and $g(\bu_k)<g(\bu_i)=f(\bv_1)\le
f(\bu_k)$ for all $s<k<i$, since $\bv_1\prec\pi(\bu_k)$ on
$\pi(W_{i,J})$. Hence, the segments $[\bu_l,\bu_{l+1}]$, $s\le
l<i$, lie in the domain $\{f\ge g\}$, and thus are contained in
$Z(F)$. \smallskip

\item If $U\cap
 \Sig _{\bu_k}(K)\ne\emptyset$ for some $k=1,\dots,i$ and
$K=\{1,3\}$,  then $\pi(W_{i,J})$ contains the negatively directed
ray starting at $\pi(\bu_k)$ and parallel to the $\var_1$-axis.
Hence, the value $f(\bu_k)$ is attained by at least two monomials
involving  $\var_2$, which keep their value along the negatively
directed ray starting at $\bu_k$ and parallel to the
$\var_1$-axis. As we have seen earlier $f(\bu_k)\ge g(\bu_k)$, and
thus the latter ray  lies entirely in the domain $\{f\ge g\}$.
Hence, $\overline{ \Sig _{\bu_k}(K)}\subset Z(F)$. The case of
$K=\{2,3\}$ is treated in the same way. \smallskip

\item For each $k=1,\dots,i-1$, the value of $g$ along
$\overline{ \Sig _{\bu_k}(J)}$ is attained by two monomials of
$g$, and thus,
$$\begin{array}{lll}
\overline{ \Sig _{\bu_k}(J)}\cap Z(F)
 & =& \left(\overline{ \Sig _{\bu_k}(J)}\cap\{g\ge
f\}\right)\cup\left(\overline{ \Sig _{\bu_k}(J)}\cap\{f\ge g\}\cap
Z(f)\right) \\[1mm]
& \supset& \overline{ \Sig _{\bu_k}(J)}\cap U\ .
\end{array}$$

\item Finally, the value of $g$ along $\overline{ \Sig _{\bu_i}(J)}$
is attained exactly by one monomial of $g$, and hence
$$\overline{ \Sig _{\bu_k}(J)}\cap Z(F)=\left(\overline{ \Sig _{\bu_k}(J)}\backslash\{g>f\}\right)
\cap\{f\ge g\}\cap Z(f)\ .$$
\end{enumerate}

Recall that  $U$ contains two negatively directed rays, starting
at $\bv_1$ and parallel to the $\var_1$-axis and to the
$\var_2$-axis, respectively, as addressed  in Section \ref{secn4}.
Now,  since the value $g(\bu_i)=f(\bv_1)$ is attained by four
monomials of $f$, two in $\var_1$ and two in $\var_2$ , we
conclude  that
$$\pi(\overline{ \Sig _{\bu_k}(J)}\backslash\{g>f\})=Q_{i,J}$$ (see
the definition at the beginning of Step 2).

Summarizing these conclusions, we have shown that $Z(F)\supset U$
for all suitable generators $f$ of the simple ideal of
$\pi(W_{i,J})$, and thus, due to (\ref{nne74}), that
$$\bigcap_fZ(F)\cap\overline{ \Sig _{\bu_i}(J)}=U_{i,J}\
.$$

\medskip

{\it Step 4}. In the case when $\hat V_{i,J}$ has no terminal
vertex we shall suitably modify the preceding construction of the
polynomials $f(\var_1,\var_2)$ and $g(\var_3)$, constituting $F$,
and at the end we shall append an additional simple polynomial
$F'$ meeting requirements (\ref{nne73}).

Consider the  additive tropical plane curve $\pi(U)\subset\Real^2$
and denote the minimal vertex of $\left(\pi(U)\right)_+$  by
$\bw=(w_1,w_2)$. Note that $\bw\prec\pi(\bu_k)$ for all
$k=1,\dots,m+1$, and that
\begin{equation}\pi(U)\cap\{\var_1<w_1,\ \var_2<w_2\}=
V_{i,J}\cap\{\var_1<w_1,\ \var_2<w_2\}\label{e9}\end{equation}
(this is just an open ray).

Again, using Propositions \ref{p1} and \ref{p2}, we can choose a
simple polynomial $f(\var_1,\var_2)$ defining $\pi(W_{i,J})$ and
satisfying the following conditions:
\\
 the sequence
$$f(\bw)<f(\bu_{i+1})< \cdots <f(\bu_{m+1})$$
is convex with respect to the sequence $u_{i3}< \al _{i+1,3}<
\cdots <u_{m3}<u_{m+1,3}$, and relation (\ref{e15}) holds true as
well. Then, we define  the polynomial $F$ as in formula
(\ref{e8}), where the arguments $A_k$ and $c_k$, $0\le k\le m+1$,
are determined by conditions (a) and  (b) in Step 3 and by the
following requirements:
\begin{enumerate}\item[(c')] for $i$-th monomial,
$$c_i=\frac{f(\bu_{i+1})-f(\bw_1)}{u_{i+1,3}-u_{i3}},\quad
\left(A_i\TrP \var_3^{c_i}\right)\Big|_{\bu_{i+1}}=
\left(A_{i+1}\TrP
\var_3^{c_{i+1}}\right)\Big|_{\bu_{i+1}}=f(\bu_{i+1})\ ,$$
$$\left(A_i\TrP \var_3^{c_i}\right)\Big|_{\bu_i}=f(\bw)\ ;$$
\item[(d')] for $0\le k<i$,
$$c_k=c_{k+1}-1,\quad\left(A_k\TrP
\var_3^{c_k}\right)\Big|_{\bu_{k+1}}=\left(A_{k+1}\TrP
\var_3^{c_{k+1}}\right)\Big|_{\bu_{k+1}}\ .$$
\end{enumerate} Conditions (a), (b), (c'), and  (d')  determine
uniquely the values of the arguments  $A_k$ and $c_k$, $0\le k\le
m+1$, out of the values
\mbox{$f(\bw),f(\bu_1),\dots,f(\bu_{m+1})$}.

Using the argument of Step 3, where $\bw$ plays the role of
$\bv_1$, we show that $Z(F)\supset U$. However, in the quadrant
$\overline{ \Sig _{u_i}(J)}$ we obtain
$$Z(F)\cap\overline{ \Sig _{\bu_i}(J)}=\left(\overline{ \Sig _{\bu_i}(J)}\cap\{g\ge
f\}\right)\cup\left(\overline{ \Sig _{\bu_i}(J)}\cap Z(f)\right)\
,$$ since the value $g(\bu_i)=f(\bw)$ is attained by two monomials
of $g$. Here, $\{g\ge f\}$ cuts off $\overline{ \Sig _{\bu_i}(J)}$
the quadrant $Q=\{\var_1\le w_1,\ \var_2\le w_2,\
\var_3=u_{i3}\}$, and $Z(f)$ cuts off $\overline{ \Sig
_{\bu_i}(J)}\backslash Q$ the set $(\overline{ \Sig
_{\bu_i}(J)}\backslash Q)\cap U_{i,J}$. Hence,
$$\overline{ \Sig _{\bu_i}(J)}\cap\bigcap_fZ(F)\cap
Z(F')=U_{i,J}$$ as required, where $F'(\var_1,\var_2)$ is a simple
polynomial defining the tropical set-curve $\pi(U)$.

The proof of Theorem \ref{thm:end} is completed.
\end{proof}

\section{Additive tropical subvarieties of additive tropical
varieties}\label{sec-last}

\begin{theorem}\label{nt11}
Let $P\subset\Real^n$ be an $m$-dimensional additive tropical
set-variety. Then \begin{itemize}\item the faces of the closures of
the connected components of $\Reg(P)$ define an FPC structure $ \mcP
$ on $P$; \item for any $k=0, \dots ,m-1$, the $k$-skeleton
$P^{(k)}=\bigcup_{ \sig \in \mcP ,\dim( \sig ) \le k} \sig $ and
each of its connected components is an additive tropical
set-variety.\end{itemize}
\end{theorem}

\begin{proof} Let $K_1, \dots ,K_N$ be all the connected components of
$\Reg(P)$. We extend the result of Lemma \ref{ln1} by showing that
$$ \olK_i\cap \olK_j\subset\partial\olK_i\cap\partial \olK_j,\quad \text{for every } 1\le i<j\le N\ .$$
Arguing on the contrary, in view of lemma \ref{ln1}, we assume
that $$\tau=\Inte( \olK_i)\cap\bigcup_{j\ne i}\overline
K_j\ne\emptyset,\quad\dim (\tau) =k<m-1\ ,$$
 for some $i\ne j$.  Let $\bx\in\tau$ be a
generic point. Then, there is an $(n-k)$-plane $$\Pi=\{x_{i_1}=
\dots =x_{i_k}=0\}\subset\Real^n$$ which is transverse to
$\Real\tau$. This $(n-k)$-plane is an additive tropical
set-variety, and so is $P'=P\cap(\bx+\Pi)$.

Let $K'_i$ be the germ of $\Inte( \olK_i)\cap\Pi$ at $\bx$, which
is the germ of an affine space of dimension $m-k\ge 2$, and let
$K''_i$ be the germ  of $\Pi\cap\bigcup_{j\ne i} \olK_j$ at $\bx$
such that $\dim (K''_i) =m-k\ge 2$. Furthermore, we may assume
that $K'_i$ is the whole affine $(m-k)$-space, and $K''_i$ is a
cone with vertex $\bx$.

Given a neighborhood of $\bx$, $P'$ is the union of the germ
$K'_i$ at $\bx$ and of the germ $K''_i$ at $\bx$ such that
$K'_i\cap K''_i=\bx$.  Clearly there is $\bu\in
K'_i\backslash\{\bx\}$ such that $\bx\oplus\bu\ne\bx$. Then
$\bx\not\in\bu\oplus(K'_i\cup K''_i)$, and, due to
$\bu\oplus\bu=\bu\in K'_i$, we conclude that $\bu\oplus(K'_i\cup
K''_i)\subset K'_i\backslash\{\bx\}$. By the same token, we see
that $\bx\oplus K''_i=\bx$, or, equivalently, that
$K''_i\subset\bx+\Real^n_{\le0}$, which in turn contradicts the
fact that $K''_i$ is a tropical set-variety. (Clearly, $K'_i$ does
not affect the balancing condition for cells of $K''_i$.)

\medskip

Now we show that $P^{(m-1)}$ (see Definition \ref{d34}) is an
additive tropical set-variety. Pick a set $I\subset\{1,2, \dots
,n\}$, $|I|=m+1$, and introduce the projection
$$\pi_I:\Real^n\to\Real^I=\{x_i=0,\ i\not\in I\} \ . $$ If, for
some $j=1, \dots ,N$, $\dim \big( \pi_I( \olK_j) \big)=m$, then
$\pi_I(P)$ is an $m$-dimensional tropical set-hypersurface in
$\Real^I$ (i.e., push-forward,  as defined in \cite{AR}) which is
additive, since $\pi_I$ is a semigroup homomorphism. Furthermore,
$Q_I=\pi_I^{-1}(\pi_I(P))$ is an additive tropical
set-hypersurface in $\Real^n$ (i.e., pull-back,  as defined in
\cite{AR}) which can be viewed  as the union of the
$(n-1)$-dimensional polyhedra $ \olK_i+V_I$, where $\dim \big
(\pi_I( \olK_i) \big)=m$ and $V_I=\{x_j=0,\ j\in I\}$.

Observe that the  two $(n-1)$-dimensional polyhedra $ \olK_i+V_I$
and $ \olK_j+V_I$ cannot intersect so that
$$\dim\left( \ \Inte( \olK_i+V_I)\cap\Inte( \olK_j+V_I)\ \right)=n-2\ .$$ Indeed, otherwise the
$(n-2)$-cell $( \olK_i+V_I)\cap( \olK_j+V_I)$ of $Q_I$ would be
dual to a parallelogram in the subdivision $S(f)$ of the Newton
polytope of a simple tropical polynomial $f$ with $Z(f)=Q_I$ (see
Theorem \ref{nt1}), and thus we reach a contradiction, since no
four distinct points on coordinate axes may span a parallelogram.

Consider now the intersection $R_{I,t}=P\cap(at+Q_I)$, for a
generic vector $a\in\Real^n$ and a small positive parameter  $t$.
This intersection is transversal and makes $R_{I,t}$ an additive
tropical set-variety of dimension $m-1$. The top dimensional cells
of $R_{I,t}$ appear as the intersections $ \olK_i\cap(at+\overline
K_j+V_I)$, $i\ne j$, where $ \olK_i$ and $ \olK_j+V_I$ are
transverse and intersect on their interiors. As $t\to0$, such an
intersection either contracts or converges to an
$(m-1)$-dimensional polyhedron in $P$. Two situations may occur:
\begin{itemize}\item If $ \olK_i\cap \olK_j$ is a common
$(m-1)$-dimensional face $ \sig $, then $\olK_i\cap(at+
\olK_j+V_I)$ converges to $ \sig $ as $t\to0$.
\item If $\dim( \olK_i\cap \olK_j)<m-1$, but $\dim(\overline
K_i\cap(at+ \olK_j+V_I))=m-1$, we necessarily obtain that $
\olK_i+V_I$ and $ \olK_j+V_I$ intersect transversally along their
interior, which is not possible as observed in the preceding
paragraph.
\end{itemize} Thus, we get that the limit $R_I$ of $R_{I,t}$ as
$t\to0$ is an additive tropical set-variety which is the union of
some $(m-1)$-faces of $ \olK_i$, $i=1, \dots ,N$. Noticing that
$R=\bigcup_IR_I=P^{(m-1)}$, where $I$ runs over all
$(m+1)$-subsets of $\{1, \dots ,n\}$, we prove that $P^{(m-1)}$ is
an additive tropical set-variety.

\medskip

The rest of the proof goes by induction on descending dimension.
\end{proof}

\begin{corollary}
Given an additive tropical set-variety, the connected components
of its skeletons are contractible.
\end{corollary}

We finish with describing additive tropical set-varieties having a
disconnected skeleton of a positive dimension:

\begin{lemma}
If a  connected additive $m$-dimensional tropical set-variety $P$ in
$\Real^n$ has a disconnected skeleton $P^{(d)}$, $0<d<m$, then there
are additive transversal linear subspaces $U,V\subset\Real^n$ of
dimension $d$ and $n-d$ respectively, such that $P=U+Q$, with
$Q=P\cap V$ being an additive $(m-d)$-dimensional tropical
set-variety.
\end{lemma}

\begin{proof}
We use the following fact which we leave to the reader as an
elementary exercise: if the union of the $d$-dimensional faces of
an $m$-dimensional convex polyhedron $ \sig $ in $\Real^n$
($n>m>d$) is not connected, then there are transversal linear
subspaces $U',V'\subset\Real^n$ of dimension $d$ and $n-d$
respectively, such that $ \sig =U'+\tau$, $\tau= \sig \cap V'$
being an $(m-d)$-dimensional convex polyhedron.

We can assume that $P^{(d+1)}$ is connected. Let $ \sig $ be a
$(d+1)$-cell of $P^{(d+1)}$ joining two connected components of
$X^{(d)}$. As noticed above, $ \sig =U+\tau$ with $\tau$ a
segment. Then $\partial \sig $ is the union of two affine spaces
$a+U$, $b+U$, $a,b\in\Real^n$. By Theorem \ref{nt11}, $P^{(d)}$ is
a tropical set-variety, and hence, due to the balancing condition,
$a+U$, $b+U$ must be separate connected components of $P^{(d)}$.
This immediately implies that $P^{(d)}$ is the union of several
affine spaces $a+U$. Since each of them is additive (Theorem
\ref{nt11}), $U$ is additive. Again the above observation on
polyhedra with a disconnected $d$-skeleton shows that each cell $
\sig $ of $P$ of dimension $>d$ is represented as $U+\tau$, and
hence $P=U+Q$, where $Q$ we can obtain as intersection of $P$ with
a transversal to $U$ additive $(n-d)$-dimensional linear subspace
$V\subset\Real^n$.
\end{proof}



\end{document}